%
%
%
%
\input amstex
\input amssym
\documentstyle{amsppt}
\def\emptyset{\varnothing}
\def\e{\epsilon}
\def\a{\alpha}
\def\b{\beta}

\def\d{\delta}
\def\G{\Gamma}
\def\g{\gamma}
\def\k{\kappa}

\def\N{\Bbb N}
\def\l{\lambda}

\def\s{\sigma}

\def\x{\times}
\def \R{\Bbb R}
\def \Z{\Bbb Z}

\def\o{\overline}
\def\f{\flushpar}
\def\u{\underline}

\def\om{\omega}
\def\Om{\Omega}
\def\B{\Cal B}

\def\({\biggl(}
\def\){\biggr)}

\def\<{\bold\langle}
\def\>{\bold\rangle}

\def \r{\Cal R}
\def\Dom{\text{\rm Dom}\,}
\def\Im{\text{\rm Im}\,}\document\topmatter
\def\wh{\widehat}
\def\textrm{\text}
\def\em{\it}
\def\mathbb{\Bbb}
\title exchangeable measures for subshifts\endtitle
\author J. Aaronson,\   H. Nakada \ and\  O. Sarig\endauthor
 \address[Jon. Aaronson]{\ \ School of Math. Sciences, Tel Aviv University,
69978 Tel Aviv, Israel.}
\endaddress
\email{aaro\@tau.ac.il}\endemail\address[Hitoshi Nakada]{\ \ Dept.
of Math., Keio University,Hiyoshi 3-14-1 Kohoku,
 Yokohama 223, Japan}\endaddress
\email{nakada\@math.keio.ac.jp}\endemail\address[Omri Sarig]{\ \
Dept. of Math., Penn State University, University Park, PA 16802,
U.S.A }\endaddress \email{sarig\@math.psu.edu}\endemail\abstract
Let $\Om$ be a Borel subset of $S^\Bbb N$ where $S$ is countable. A measure is called exchangeable on $\Om$, if it is supported on $\Om$ and is invariant under every Borel automorphism of $\Om$ which permutes at most finitely many coordinates.
De-Finetti's theorem characterizes these measures when $\Om=S^\Bbb N$. We apply the ergodic theory of equivalence relations to study the case $\Om\neq S^\Bbb N$, and  obtain versions of this theorem when $\Om$ is a countable state Markov shift, and when $\Om$ is the collection of beta expansions of real numbers in $[0,1]$ (a non-Markovian constraint).
\endabstract \thanks\copyright 2004.\endthanks
\endtopmatter\rightheadtext{exchangeability}\heading {\S0 Introduction}\endheading

\noindent
{\bf Exchangeability.}
De-Finetti's theorem says that if a stochastic process $\{X_n\}_{n\geq 1}$ is {\it exchangeable}, i.e. all finite permutations $\{X_{\pi(n)}\}$ of $\{X_n\}_{n\geq 1}$ are distributed like $\{X_n\}$, then it is distributed as a mixture of i.i.d. distributions.

Here is a seemingly stronger, but equivalent formulation: Let
$\frak {K}$ be the collection of all bi-measurable bijections
 $\kappa:A\to B$ ($A,B\subseteq S^\N$ Borel) for which for every $x$
$\kappa(x)$ is some finite permutation\footnote{A permutation $\pi$ is called {\it finite} if its support $\{s: \pi(s)\neq s\}$ is finite.} of $x$; then any Borel probability measure $m$ on $\Om:=S^\N$ such that $m\circ\kappa|_{\Dom(\kappa)}=m|_{\Dom(\kappa)}$ for all $\kappa\in\frak{K}$ is  an average of Bernoulli measures.

De-Finetti's theorem is instrumental in statistical modeling of sequential sampling, because it determines the form of joint distributions whenever the sampling order is unimportant. But sometimes the sampling order is subject
to  non-permutation invariant deterministic constraints. In these cases the joint distribution cannot be assumed to be exchangeable. Nevertheless, one can still ask for `the most exchangeable' compatible distributions.

There are various ways to formalize this. In this paper we use the following: Let
$\Om$ be a Borel subset of $S^\Bbb N$ (thought of as the space of realizations of $\{X_n\}_{n\geq 1}$ subject to a collection of deterministic constraints), and set $\frak{K}(\Om):=\{\kappa\in\frak{K}:\Dom(\kappa), \Im(\kappa)\subseteq\Om\}$. A Borel measure $m$ on $\Om$  is called {\it exchangeable} on $\Om$ if $m\circ\kappa|_{\Dom(\kappa)}=m|_{\Dom(\kappa)}$ for all $\kappa\in\frak{K}(\Om)$. When
$\Om=S^\N$, this reduces to the usual notion of exchangeability.

This definition of exchangeability is the one used by Petersen \& Schmidt in the context of finite state Markov shifts [Pe-S], but is not equivalent to the definition of `partial exchangeability' introduced by Diaconis \& Freedman in the context of topological Markov shifts [D-F]. (Topological Markov shifts are sample spaces of Markov chains, see \S 3 below.)

The shift invariant exchangeable measures for two sided finite state topological Markov shifts were determined by Petersen \& Schmidt [Pe-S].  The exchangeable measures for a one sided finite state topological Markov shift were determined by Aaronson, Nakada, Solomyak \& Sarig in [ANSS1]. The partially exchangeable measures for countable state Markov topological shifts were determined by Diaconis \& Freedman [D-F].

\medskip
\noindent
{\bf Aim.\/}
This paper describes  exchangeable
measures on $\Om\subset S^\N$, $|S|\leq \aleph_0$, in the following cases:
\roster
\item {\it Markov Constraints:\/} $\Om$ is a one-sided countable state topological Markov shift (i.e. the sample space of a countable state Markov chain, see \S3);
\item {\it A Non-Markov constraint:\/}  $\Om$ is a $\beta$--shift (i.e. the collection of all (greedy) $\beta$-expansions of real $\theta\in [0,1]$, where $\beta>1$, see \S7). The motivation for studying $\beta$--shifts comes from number theory (see [Re], [Pa], [Schw] and references therein).
\endroster
Our results apply to probability measures as well as to locally finite
 infinite measures (see below). Such measures appear naturally in our
context, because the state space $S$ is infinite.

\medskip
We give a brief outline of our approach.

\subheading{Equivalence relations}
The sets $\Om$ considered above are shift invariant: $T(\Om)=\Om$, where $T$ is the left shift map $T(x_1,x_2,\ldots )=(x_2,x_3,\ldots)$.
We will use the language of equivalence relations reviewed below to formulate the exchangeability property in terms of some natural equivalence relations associated with certain skew--products over $T$ (see [Pe-S], [ANSS] and below). This will allow us to bring in some tools from ergodic theory and thus bypass some of the combinatorial complications a direct approach would have encountered.

Let $(X,\B(X))$ be a standard measurable space.
An {\it equivalence relation} on
$X$ is a set $\Cal R\subseteq X\x X$ such that the relation $x\sim y\Leftrightarrow (x,y)\in\Cal R$ is an equivalence relation. An equivalence relation
is called {\it Borel} if $\Cal R\in\B(X)\otimes\B(X)$,  and is called {\it countable} if all its equivalence classes $\Cal R_x:=\{y:(y,x)\in\Cal R\}$ $(x\in X)$ are countable.
Take as an example a countable discrete set $S$, a Borel set $\Om\subseteq S^\N$ with the relative product topology, and  $\B(\Om)$ the Borel $\s$--algebra. The following set is a countable Borel equivalence relation on $\Om$:
$$
\Cal E(\Om):=\{(x,y)\in \Om\x\Om: x,y\text{ differ by a finite permutation}\}.
$$
We later refer to $\Cal E(\Om)$ as the  {\it exchangeable relation } of $\Om$.

A bi-measurable bijection $\kappa$ defined on some $A\in {\Cal B} $
with image $B\in{\Cal B} $ is an $\Cal
R$-{\it holonomy } if $ (x,\kappa(x))\in \Cal R$ for any $ x\in A$.
 We write in this case $A\overset{\Cal R}\to\rightarrow B$.

A function $F:X\to\R$ is called $\r$--invariant, if it is invariant under all $\r$--holonomies. A measure $m$ is called $\r$--ergodic, if all measurable $\r$--invariant functions are equal a.e. to a constant. Every $\r$--invariant measure can be decomposed into ergodic components, see \S1 for details.

A measure on $(X,\B(X))$ is called {\it $\Cal R$--invariant}, if $m\circ\kappa|_{\Dom(\kappa)}=m|_{\Dom(\kappa)}$ for all $\Cal R$--holonomies $\kappa$ (with domain $\Dom(\kappa)$).
The collection of $\Cal E(\Om)$--invariant measures is exactly the collection of exchangeable measures on $\Om$.

\subheading{Exchangeability, Skew-Products, and Tail Relations}
We represent the exchangeable relation in terms of $T:\Om\to\Om$.  Fix $a_0\in S$, consider the (additive) Abelian group
$$
\Bbb Z^{S\setminus
\{a_0\}}_0:=\{x\in\Bbb Z^{S\setminus \{a_0\}}: \text{ all but finitely many coordinates of $x$ are zero}\}
$$ equipped with the discrete topology, and  define
$
F^\natural:\Om\to\Bbb Z^{\a\setminus \{a_0\}}_0$ by
$$
F^\natural(x_0,x_1,\ldots)_a:=\d_{a,x_0}.
$$
Let $F_k^\natural:=F^\natural+F^\natural\circ T+\cdots F^\natural\circ T^{k-1}$. These count the appearances of elements of $S$ in the first $k$ symbols of $x$. It is routine to verify that the exchangeable relation is the same as
$$
\goth T(T,F^\natural):=\{(x,y)\in\Om\x\Om: \exists k\ge 0\text{ s.t. }T^kx=T^ky\text{ and }F^\natural_k(x)=F^\natural_k(y)\}.
$$
It follows that exchangeability is the same as $\goth T(T,F^\natural)$--invariance.

In order to study $\goth T(T,F^\natural)$, we represent it in terms of the tail relation of a suitable transformation. This is done using standard abstract ergodic theoretic constructions which we now review.
Let $X$ be a standard space.
\roster

\item {\it Tail relations:\/}  Let $T:X\to X$ be a measurable, locally invertible
transformation on $X$.
  The {\it tail relation} of $T$ is
$$
\goth T(T):=\{(x,y)\in X\x X:\ \exists\ k\ge 0\text{ s.t. }
T^kx=T^k y\}.
$$
Define $X_0:=\{x\in X: x\text{ is not eventually periodic}\}$.
The {\em grand tail relation} of $T$ is the equivalence relation
$$
\align
\goth G(T)&:=\{(x,y)\in X_0\x X_0: \exists k,\ell\ge 0\text{ s.t. }T^k x=T^\ell y\}.
\endalign
$$

\item {\it Orbit cocycles:\/} Let $\r$ be a countable Borel equivalence relation on $X$ and $\Bbb G$ an Abelian topological group.  A Borel function $\Psi:\Cal R\to\Bbb G$ is called an $\r$--{\it cocycle} if
$$
\Psi(x,z)=\Psi(x,y)+\Psi(y,z)\text{
whenever  }(x,y),\ (y,z)\in\r.
$$
Any $F:X\to\Bbb G$ gives rise to the following  $\goth T(T)$ and $\goth G(T)$--cocycles (which abusing notation we denote by the same symbol):
$$
\align
\widehat{F}(x,y)&:=F_k(y)-F_k(x),\ \ (x,y)\in \goth T(T)\text{ and }T^k x=T^k y\\
\widehat{F}(x,y)&:=F_k(y)-F_\ell(x),\ \ (x,y)\in \goth G(T)\text{ and }T^\ell x=T^k y, \tag{0.1}
\endalign
$$
where $F_k:=F+F\circ T+\cdots +F\circ T^{k-1}$.

\item {\it Skew--Products:\/}
The {\it skew--product relation} with {\it base} $\r$ and {\it cocycle} $\Psi$ is the following countable equivalence relation on $X\x\Bbb G$:
$$
\r_\Psi:=\bigl\{\bigl((x,t),(y,s)\bigr)\in (X\x\Bbb G)^2: (x,y)\in \r\text{ and }t-s=\Psi(x,y)\bigr\}.
$$
Note that $\goth T(T)_{\wh{F}}\equiv \goth T(T_F)$, where $T_F:X\x\Bbb G\to X\x\Bbb G$ is the skew--product {\em transformation} $T_F(x,\xi):=(Tx,\xi+F(x))$.

\item {\it Inducing:\/} Let $\r$ be as above and $E\subset X$ be some Borel
set. The {\it induced relation on }$E$ is $\r(E):=\r\cap (E\x E)$. This is a countable Borel equivalence relation on $E$.
\endroster
We now have the following identity:
$$
\Cal E(\Om)=\frak T(T,F^\natural)\cong \frak T(T)_{\widehat{F^\natural}}\cap [\Om\x\{0\}]^2=\frak T(T_{F^\natural})([\Om\x\{0\}]^2), \tag{0.2}
$$
where the isomorphism $\cong$ is $(x,y)\leftrightarrow\bigl((x,0),(y,0)\bigr)$.

It is a standard fact that if $\r$ is a countable Borel equivalence relation on $X$ and $E\subseteq X$ is  Borel, then any $\r$--ergodic invariant measure restricts to an $\r(E)$--invariant ergodic measure on $E$, and that any $\r(E)$--ergodic invariant measure arises this way (see Proposition 1.0 below).

This fact, together with identity (0.2) reduces the study of exchangeable measures on $\Om$ to the study of $\goth T(T_{F^\natural})\equiv (\goth T(T))_{\wh{F^\natural}}$--invariant measures on $\Om\x\Bbb Z^{S\setminus\{a_0\}}$.

\subheading{Conformal measures and the Maharam construction}
The previous discussion shows that the exchangeability problem can be reduced to the study of invariant measures for the skew--product relation $\r_\Psi$, with $\r=\goth T(T)$ and $\Psi=\wh{F^\natural}$. There is a standard construction, called the {\it Maharam construction} after [Ma], of such measures.

Let $\r$ be a countable Borel equivalence relation on $X$, $\Bbb G$ a locally compact polish Abelian group (e.g. $\Bbb Z^{S\setminus\{a_0\}}$),  $\Psi:\r\to\Bbb G$ an orbit cocycle, and $H:\Bbb G\to\R$ a continuous homomorphism. A measure $\mu$ on $X$ is
{\it $\r$--non-singular}, if every $\r$--holonomy $\kappa:A\to B$ is non-singular, i.e., $m\circ\kappa|_A\sim m|_A$.
 An $\r$--non-singular measure $\mu$ is called $(e^{H\circ\Psi},\Cal
R)$-{\it conformal} if
$\tfrac{d\mu\circ\kappa}{d\mu}(x)=e^{H\circ\Psi(x,\kappa x)}$ a.e. on  $\Dom \k$ for all holonomies $\kappa$.

The {\em Maharam measure} corresponding to a $(e^{H\circ\Psi},\r)$--conformal measure $\mu$ on $X$ is the following measure on $X\x\Bbb G$:
$$
dm(x,y):=e^{-H(y)}d\mu(x)dm_\Bbb G(y),
$$
where $m_\Bbb G$ is a Haar measure for $\Bbb G$. It is straightforward to check that such measures are $\r_\Psi$--invariant. Ergodicity is {\em not} guaranteed even when $\mu$ is $\r$--ergodic.

It is instructive to interpret the Maharam measures in the special case $X=\Om$, $\Bbb G=\mathbb Z^{S\setminus\{a_0\}}$,  $\r=\goth T(T)$, and $\Psi=\wh{F^\natural}$. In this case any Maharam  measure on $\goth T(T)_{\wh{F^\natural}}$ restricts to the measure $\mu\x\d_{0}$ on $\goth T(T)_{\wh{F^\natural}}\cap (\Om\x\{0\})^2$. The isomorphism $\cong$ carries this measure to  $\mu$, and identity (0.2) shows that this measure must be exchangeable.

Thus every $(e^{H\circ F^\natural},\goth T(T))$--conformal measure is exchangeable.
In this paper we study the other direction: starting with an arbitrary exchangeable measure, we ask to what extent can it be constructed from conformal measures. We do this for two particular choices of $\Om$: countable Markov shifts, and sets of $\b$--expansions.

\subheading{Programme}
For the sets $\Om$ described above we shall do the following:
\roster
\item Establish the existence of conformal measures, and identify them;
\item Characterize their $\goth T(T,F^\natural)$--ergodicity;
\item Show, or find sufficient conditions for $\goth T(T,F^\natural)$--invariant ergodic measures to be ergodic conformal measures when restricted to their support.
\endroster

We compare the Markovian and non-Markovian situations.
The study of conformal measures in the non-Markovian case ($\b$--expansions) requires different tools than in the Markovian case (countable Markov shifts).
There is however a common thread in part (3).

Rather than showing that $\goth T(T,F^\natural)$--ergodic invariant locally finite measures are conformal, we show that $\goth T(T)_{\wh{F^\natural}}$--ergodic invariant locally finite measures $m$ are Maharam. The key step is to show quasi--invariance under all transformations of the form $Q_a(x,\xi):=(x,\xi+a)$, as this implies the Maharam form for abstract reasons (see the proof of Theorem 5.0 or 9.0). To do this we construct explicit approximations to $Q_a$ by $\r_\Psi$--holonomies.

The construction of the approximating holonomies  depends, of course, on the structure of the $\Om$ in question.

\medskip

Various parts of this programme make sense for more general cocycles than $F^\natural$, and a larger class of subshifts than those considered here. It is of some interest to identify particular properties which are sufficient for our argument to work. We therefore carry out parts of this programme in greater generality than needed for the  exchangeability per-se.

The paper is divided into three parts. In the first, we collect some terminology, notation, and facts from the ergodic theory of equivalence relations that will be needed in the sequel. We also solve the exchangeability problem for locally finite {\em infinite} measures on full--shifts. In the second and third parts we treat, respectively, exchangeable measures for countable Markov shifts, and for $\b$--expansions.

\heading{Part I: Generalities}
\endheading

\heading{\S1 More on equivalence relations and conformal measures}
\endheading

\subheading{Some notions of finiteness for measures}
Let $(X,\Cal B(X))$ be a standard measurable space.
The
collections  of probabilities, and $\s$-finite measures on a
standard measurable space $X$ are denoted by $\Cal P(X)$ and $\goth
M(X)$ respectively. If $\a\subset\B(X)$ is a countable
partition, we set $\goth M_\a(X):=\{\mu\in\goth M(X):\
\mu(A)<\infty\ \forall\ A\in \a\}$ and call these measures {\it
$\a-\s$-finite}.

If $X$ is equipped with a topology (generating
its measurable structure), we call $\mu\in\goth M(X)$ {\it
locally finite} (on $X$) if there is a countable cover  of $X$
by open sets, each with  finite $\mu$-measure, and
 {\it
topologically $\s$-finite} if it is locally finite on some
Borel subset of full  $\mu$-measure.

\subheading{The Feldman--Moore Theorem}
The definition of ergodicity and invariance under an equivalence relation includes quantification over all holonomies. This is unnecessary.

A collection $\goth C$ of $\Cal R$- holonomies {\it generates
$\Cal R$} if for each $(x,y)\in\Cal R,\ \exists\  \Phi\in\goth C$
such that $x\in\Dom\Phi,\ y=\Phi(x)$. A measure is $\r$--invariant iff it is invariant w.r.t. to a generating collection of holonomies. A function is $\r$--invariant iff it is invariant w.r.t a generating collection of holonomies.

Feldman \& Moore proved in \cite{F-M} that any countable Borel equivalence relation is generated by a {\em countable} group of globally defined holonomies: $\exists\Gamma$ a group of Borel automorphisms s.t. $\r=\r_\Gamma:=\{(x,g(x)):x\in X, g\in\Gamma\}$.

\subheading{The ergodic decomposition}
For a countable Borel equivalence relation,
a set $A\subset X$ is  $\r$-invariant iff  $x\in A\Longrightarrow \r_x\subseteq A$ (where $\r_x$ is the equivalence class of $x$).
The collection of $\r$--invariant sets forms a $\s$--algebra, which is
 denoted by $\frak{I}(\r)$.
A measure $\mu$ is  {$\r$--ergodic} iff every $A\in\frak{I}(\r)$ is equal to $\emptyset$ or $X$ up to a $\mu$--null set. We write in this case
$\frak{I}(\r)\overset{\mu}\to{=}\{\emptyset,X\}$.

\par A {\it mixture of measures} is a measure $\mu\in\goth M(X)$ of form
$$\mu(B)=\int_\Om\mu_\om(B)d\nu(\om)\ \forall\ B\in\B(X)$$
where $(\Om,\Cal F,\nu)$ is a $\s$-finite measure space,
$\mu_\om\in\goth M(X)\ \forall\ \om\in\Om$ and $\om\mapsto\mu_\om$
is measurable ($\Om\to\goth M(X)$) in the sense that
$\om\mapsto\mu_\om(A)$ is  measurable $\forall\ A\in\B(X)$.

 The mixture is called {\it finite} or
{\it infinite} according to whether  $\nu(\Om)<\infty$ or
$\nu(\Om)=\infty$, and an {\it average} if $\nu(\Om)=1$. Note that
an average of finite measures could also be an infinite mixture of
probabilities (see the examples after proposition 2.2). The measures
$\{\mu_\om:\ \om\in\Om\}$ are called {\it components} of (the
decomposition of) $\mu$ (note that  components are only specified
up to $\nu$-measure zero).

If $\r$ is a countable Borel equivalence relation, then any
 $\r$-invariant $\mu\in\goth M(X)$  is an average of $\r$-invariant,
ergodic, $\s$-finite measures \cite{G-S}.

This reduces the problem of classifying $\Cal E(\Om)$--invariant measures to that of classifying  {\it ergodic} $\Cal E(\Om)$--invariant measures.

\subheading{Inducing and invariant measures}
The following result is a useful property of ergodic measures:
\proclaim{Proposition 1.0}
If $A\in\B$ and $\mu\in\goth M(A)$ is $\r(A)\equiv \r\cap (A\x
A)$-invariant and ergodic, then there is a unique $\o\mu\in\goth
M(X)$, $\r$-invariant and ergodic, such that
$\o\mu|_A\equiv\mu$.\endproclaim

\medskip
\noindent
{\it Proof.\/}
Let $\r=\r_\G$ where $\Gamma$ is a countable group of automorphisms which generates $\r$ as in [F-M]. There are $\g_n\in\G\ \ (n\ge 1)$ and
$A_n\in\B(A)$ so that $\bigcup_{\g\in\G}\g
A=\biguplus_{n\geq 1}\g_nA_n$. The required measure is
$\o\mu(B):=\sum_{n=1}^\infty\mu(
\g_n^{-1}(B\cap\g_n
A_n)).
$
\hfill\qed

\subheading{Conformal measures for tail relations}
The conformality property can be significantly simplified when the underlying equivalence relation is a tail relation. Call $\mu\in\goth M(X)$ $(e^F,T)$--{\em conformal} if $\mu$ is $T$--nonsingular and $\frac{d\mu\circ T}{d\mu}=e^F$. The following proposition relates this notion to conformality w.r.t $\goth T(T)$ and $\goth G(T)$:
\proclaim{Proposition 1.1}Let $(X,T)$ be a measurable, locally invertible
transformation of a  measurable space, let $F:X\to\Bbb R$  be
measurable and suppose $\wh{F}$ is as in (0.1).
\roster
\item  $\mu\in\goth
M(X)$ is $(e^{\widehat F},\goth G(T))$-{ conformal}  iff it is $T$--nonsingular, its
support is $\goth G(T)$-invariant and $\tfrac{d\mu\circ
T}{d\mu}=e^F$;
\item If  $\mu\in\goth M(X)$ is $\goth G(T)$--nonsingular, and $\frac{d\mu\circ T}{d\mu}=ce^{F}$ for some constant $c>0$, then  $\mu$ is $(e^{\widehat F},\goth
T(T))$-{ conformal}.
\item If
$\mu\in\goth M(X)$ is $(e^{\widehat F},\goth
T(T))$-{ conformal}, $\goth T(T)$-ergodic, and $\goth G(T)$-nonsingular, then $\exists c>0$ constant such that
$\tfrac{d\mu\circ T}{d\mu}=ce^F\ \mod\/
\mu$.
\endroster
\endproclaim
\demo{Proof}
The first two statements are established by direct calculation. We prove the third.
 By assumption, for any $\goth T(T)$-holonomy $K$,
 $\tfrac{d\mu\circ K}{d\mu}(x)=e^{\widehat F(x,Kx)}$. But also,
$\tfrac{d\mu\circ K}{d\mu}(x)=e^{\widehat{\log
T^\prime}(x,Kx))}$, where $T':=\tfrac{d\mu\circ T}{d\mu}$. It follows that $\widehat
F\equiv\widehat{\log T^\prime}$ on $\goth T(T)$.  \par We claim
that $h:=F-\log T^\prime$ is $\goth T(T)$-invariant, whence
constant $\mu$-a.e. To see this, note first that $\widehat
h=\widehat F-\widehat{\log T^\prime}\equiv 0$ on $\goth T(T)$. For
$(x,y)\in\goth T(T)$, we have that $h_n(x)=h_n(y)$ whenever
$T^nx=T^ny$. If $T^nx=T^ny$, then also $T^{n+1}x=T^{n+1}y$ and
$h(x)=h_{n+1}(x)-h_n(Tx)=h_{n+1}(y)-h_n(Ty)=h(y).$
Thus $\tfrac{d\mu\circ T}{d\mu}=ce^{F}$ for some
$c>0$.\hfill\qed\enddemo
\noindent
Finally, here is a straightforward generalization of a calculation done in \S 0. Define
$$
\align
\goth T(T,F)&:=\{(x,y)\in X^2: \exists n\geq 0\text{ s.t. }T^n x=T^n y\textrm{, }F_n(x)=F_n(y)\};\\
\goth G(T,F)&:=\{(x,y)\in {X_0}^2: \exists k,\ell\geq 0\text{ s.t. }T^k x=T^\ell y\textrm{, }F_k(x)=F_\ell(y)\}. \tag{1.0}
\endalign
$$
\proclaim{Proposition 1.2} Let $\Bbb G$ be an
Abelian topological group and $F:X\to\Bbb G$ measurable. If
 $\mu\in\goth M(X)$ is $\goth G(T)$-non-singular and $\tfrac{d\mu\circ T}{d\mu}=ce^{H\circ F}\ \mod\/
\mu$ where $H:\Bbb G\to\Bbb R$ is a homomorphism and $c>0$, then
$\mu$ is $\goth T(T,F)$-invariant.
\endproclaim
\noindent
The proof of proposition 1.2 is immediate from the definitions.
As we shall see, much of this paper is boils down to its converses.

\subheading{The Glimm-Effros Theorem}
Consider as an example the one--sided two shift, i.e. the map
 $T$ on
$\Om:=\{0,1\}^\Bbb N$ defined by $T(x_1,x_2,\ldots )=(x_2,x_3,\ldots)$. As is well known, there exists a unique $\frak T(T)$--invariant probability measure: the $(\frac{1}{2}, \frac{1}{2})$--Bernoulli measure. Nevertheless, there are uncountably many non-atomic mutually singular $\frak T(T)$--ergodic and invariant $\s$-finite measures:  Pick $K\subset\N$ with infinite complement
 and consider the
probability measure $\mu_K$ on $\Om_K:=\{x\in\{0,1\}^\N: i\in K\Rightarrow x_i=0\}$ obtained from the $(\frac{1}{2}, \frac{1}{2})$--Bernoulli measure after the identification $\Om_K\simeq \{0,1\}^{\Bbb N\setminus K}$. Now extend it to a $\s$--finite ergodic invariant measure on $\{0,1\}^\N$ using Proposition 1.0.

The same phenomena occurs for a general countable Borel equivalence relation, as soon as it admits a non-atomic ergodic invariant measure, as explained below.  Recall that any countable Borel equivalence relation  is of the form $\r_\Gamma:=\{(x,g(x)):x\in X, g\in\Gamma\}$ for some countable group of Borel automorphisms. We have:
\proclaim{Glimm-Effros theorem\ \ \cite{E},\cite{Gl}}
Let $\G$ be a countable group of Borel automorphisms of the
standard Borel space $X$, then either (i) $X=\biguplus_{\g\in\G}\g
A$ for some $A\in\B$, or (ii) $\exists\ A\in\B$ so that
$A\cong\{0,1\}^\Bbb N$ and $\r_\G\cap A\x A\cong\goth T(T)$ where
$T$ is the shift on $\{0,1\}^\Bbb N$.\endproclaim
\noindent
If $\r$ has at least one non-atomic ergodic non-singular measure, then case (i) cannot occur, and case (ii) must hold. But in this case there are uncountably many mutually singular non-atomic $\s$--finite ergodic invariant measures,
because of the construction sketched above (see also \cite{Schm}).

 \heading{\S2 Fibred systems
and exchangeability}\endheading \subheading{Fibred systems} It is
convenient to work with measurable fibred systems, because this setting allows to describe countable Markov shifts, $\b$--expansions, and many other
subshifts  easily.

A (measurable)
{\it fibred system} is a triple $(X,T,\a)$ where $X$ is a standard
measurable space, $T:X\to X$ is a measurable transformation  and
$\a\subset\B(X)$ is a finite or countable partition  such that:
\roster
\item $\bigvee_{i=0}^\infty T^{-i}\a$ generates $\B$;
\item for every $A\in \a$, $T|_A:A\to TA$ is bi-measurable,
invertible.
\endroster

\par If $(X,T,\a)$ is a fibred system, then for each $k\in\Bbb N$,
$(X,T,\a_k)$ and $(X,T^k,\a_k)$ are also fibred systems, where
$$\a_k:=\bigvee_{j=0}^{k-1}T^{-j}\a=\{\bigcap_{j=0}^{k-1}T^{-j}a_j,\ a_0,\dots,a_{k-1}\in\a\}.$$
The elements of $\a_k$ are called {\it cylinders} (of length $k$)
and are denoted
$\bigcap_{j=0}^{k-1}T^{-j}a_j=:[a_0,a_1,\dots,a_{k-1}]$.\par If
$S$ is countable and $X\subset S^\Bbb N$ is a subshift, then
$(X,T,\a)$ is a fibred system where $T$ is the shift and
$\a:=\{a_s:\ s\in S\},\ a_s:=\{x\in X:\ x_1=s\}$. Here we abuse
notation and denote $a_s=[s]$.\par Let $Y$ be a set. The {\it
memory } of a function $f:X\to Y$ is the minimum $k\in\Bbb N\cup\{0\}$ so
that $f$ is   {\it $\a_{k+1}$-measurable} (i.e. constant on each
$a\in\a_{k+1}$); A function is said to have {\it infinite memory} if it
is not $\a_k$-measurable for any $k\ge 1$. \par Let $(Y,d)$ be a
metric space. For $N\ge 1$, we say that $\phi:X\to Y$ is
$(\a,N)$--H\"older continuous on $A\subset X$ if $\exists\rho\in
(0,1), M>0$ such that $x,y\in A\cap\a_k\Rightarrow
d(\phi(x),\phi(y))\leq M\rho^k$ for all $k\geq N$.
\par We call the fibred system $(X,T,\a)$ {\it full} if $TA=X\
\forall\ A\in\a$ and we call a probability $m\in\Cal P(X)$  a
$(X,T,\a)$-{\it product measure} if
$$m\bigg(\bigcap_{j=0}^nT^{-j}A_j\bigg)=\prod_{j=0}^nm(A_j)\ \forall\ n\ge
 0,\ A_1,\dots,A_n\in\a.$$

\par Given a fibred system $(X,T,\a)$, we  define
$\pi:X\to\a^\Bbb N$ by $T^{n-1}x\in\pi(x)_n\in\a$ and consider the
subshift  $\Sigma(X,T,\a):=\o{\pi(X)}\subset\a^\Bbb N$.
 The map $\pi:X\to\a^\Bbb N$ is injective, and we denote $x=(a_0,a_1,\dots)$
where $\pi(x)=(a_0,a_1,\dots)$. Eventually periodic points of form
$p=(a_0,a_1,\dots,a_n,b_1,\dots,b_k,b_1,\dots,b_k,\dots)$ are
denoted $p=(a_0,a_1,\dots,a_n,\o{b_1,\dots,b_k})$.

\subheading{The exchangeable relation of a fibred system} The exchangeable equivalence relation of a
fibred system $(X,T,\a)$ is
$$\Cal E(X,T,\a):=(\pi\x\pi)^{-1}\Cal E(\Sigma(X,T,\a)).$$
As before  $\Cal E(X,T,\a)=\goth T(T,F^\natural)$, where
$$
F^\natural=F^{\natural,a_0}:X\to\Bbb Z^{\a\setminus \{a_0\}}_0,\ F^\natural(x)_a:=1_a(x)\ \ \ (a\in\a\setminus a_0),
$$
$a_0\in\a$ is fixed, and  $\Bbb Z^{\a\setminus
\{a_0\}}_0:=\{x\in\Bbb Z^{\a\setminus \{a_0\}}:\ \#\{s\in \a:\
x_s\ne 0\}<\infty\}$ equipped with the discrete topology and vector addition.

We note for future reference that   $\{F^\natural(x)-F^\natural(y):\ x,y\in X\}$
generates $\Bbb Z^{\a\setminus \{a_0\}}_0$ (this would not have been true had we worked with $\Z_0^S$). We also note that $F^\natural$ is $\a$--measurable.

\subheading{Conformal measures for fibred systems}
Recall that our plan is to relate exchangeable measures to $(\goth T(T), e^{H\circ F^\natural})$--conformal measures for some continuous homomorphism $H:\Z_0^{S\setminus\{a_0\}}\to\R$.
The following proposition relates conformality with the property of having $\a$--measurable derivative (compare with  Proposition 1.1 above):
\proclaim{Proposition
2.0}
 Let $\mu\in\goth M_\a(X)$ be $\goth G(T)$-non-singular, then
\roster
\item $\tfrac{d\mu\circ T}{d\mu}$ is $\a$-measurable
$\Leftrightarrow$ $\tfrac{d\mu\circ T}{d\mu}=ce^{H\circ
F^\natural}$ for some $c>0$ and some homomorphism $H:\Bbb
Z^{\a\setminus\{a_0\}}_0\to\Bbb R$; and
\item   in this case,
$\mu$ is  $\Cal E(X,T,\a)$-invariant.
\endroster
\endproclaim
\medskip
\noindent{\it Proof.\/} The $(\Leftarrow)$ implication in (1) is clear. We prove ($\Rightarrow$):\ Let $c$ be the value of $\tfrac{d\mu\circ T}{d\mu}$ on $a_0$, and
define a homomorphism $H:\Bbb Z^{\a\setminus\{a_0\}}_0\to\Bbb R$
by $H(e_a):=\log\tfrac{d\mu\circ T}{d\mu}|_a-\log c$ (this is a constant) where for every
$a\in\a,\ (e_a)_b=\d_{a,b}\ \ (b\in\a)$. We have $\tfrac{d\mu\circ
T}{d\mu}=ce^{H\circ F^\natural}$. Part (2) follows from proposition
1.2.  \hfill\qed

\subheading{Recurrence and the de Finetti-Hewitt-Savage theorem}
The de Finetti-Hewitt-Savage theorem states:
 \proclaim{Theorem 2.1}
 Let $(X,T,\a)$ be a full fibred system.
\roster
\item $\Cal E(X,T,\a)$ is ergodic with respect to any
 $(X,T,\a)$-{product measure}.
\item If $\mu\in\Cal P(X)$ is $\Cal E(X,T,\a)$-invariant, then
 $\mu$ is an average of $(X,T,\a)$-product probability measures.
\item
 If $\mu\in\Cal P(X)$ is $\Cal E(X,T,\a)$-invariant,  then
 $\goth{I}(\Cal E(X,T,\a))\overset{\mu}\to{=}\goth{I}(\goth T(T))$.
\endroster\endproclaim
\demo{Proof} See [H-S] (also [Me] and [D-F]).\hfill\qed\enddemo
\noindent
We'll need to consider extensions of  part (2) of the theorem for
$\s$-finite measures.

For this purpose we need the following definitions. Suppose $(X,T,\a)$ is a fibred system, and define for every $a\in\a$.
$$
N_a:=\sum_{n=0}^\infty 1_{a}\circ T^n
$$
A measure $\mu\in\goth M(X)$  is called {\it recurrent}  (w.r.t. to $(X,T,\a)$) if  $N_a\in \{0,\infty\}$ $\mu$-a.e. for every $a\in\a$.

\proclaim{Proposition 2.2} Let $(X,T,\a)$ be a full fibred system, and suppose $\mu\in\goth M(X)$ is  $\Cal E(X,T,\a)$-invariant.
\roster
\item  If $\mu$ is topologically $\s$-finite and recurrent, then
 $\mu$ is a  mixture of $(X,T,\a)$-product probability measures.
\item
 If $\mu$ is locally finite, then $\mu$ is
 recurrent, whence by (1) a mixture of $(X,T,\a)$-product probability measures.
\endroster
 \endproclaim
\medskip
\noindent
{\it Proof.\/} Suppose $\mu$ is topologically $\s$-finite and recurrent $\Cal E$--invariant measure, where $\Cal E:=\Cal E(X,T,\a)$. For every $x\in X$, define $S_x:=\{x_n:\ n\ge 1\}$. By
recurrence,  $\forall\ s\in
S_x,\  N_s(x)=\infty$ for $\mu$-a.e. $x\in X$.

\par First, we consider possible atoms of $\mu$. Suppose
$x=(s_1,s_2,\dots)\in X$ is an atom. We claim that
$s_n=s_1\ \forall\ n\ge 1$. Otherwise, fix (using topological $\s$--finiteness and $\mu\{x\}\neq 0$) some $N$ such that $[s_1,\ldots,s_N]$ has finite measure. Since $\mu$ is recurrent $|S_{T^N x}|=|S_x|>1$, and consequently $(s_k)_{k\geq N}$ has an infinite number of finite permutations $\{y(k)\}_{k\geq 1}$. Using the assumption that $(X,T,\a)$ is full, we see that
$z(k):=(s_1,\ldots,s_N,y(k))\in X\cap [s_1,\ldots,s_N]$. By exchangeability, $\mu\{z(k)\}=\mu\{x\}\neq 0$. But this implies that $\mu[s_1,\ldots,s_N]=\infty$, a contradiction.

The conclusion is that any atom of $\mu$ must be of the form $x=(s,s,\ldots)$. If $x=(s,s,\dots)$ then $\d_x$ is an atomic
$(X,T,\a)$-product measure. Thus $\nu$ is a mixture of atomic
$(X,T,\a)$-product measures and a non-atomic, topologically
$\s$-finite, recurrent and
 $\Cal E$-invariant measure. This allows us to assume (as we do henceforth) that $\mu$ is non-atomic.

   By the ergodic decomposition (see \cite{G-S}), $\mu$ is an average of $\Cal
E$-invariant, ergodic components, each of which is non-atomic,
recurrent and topologically $\s$-finite. We claim each  component
is a multiple of a $(X,T,\a)$-product measure.

Accordingly, assume that $\mu\in\goth M(X)$ is topologically $\s$-finite,  recurrent and
 $\Cal E$-invariant and ergodic. For each
 $s\in S,\ \{x\in X:\ s\in S_x\}\in\goth{I}(\Cal E)$, whence
 $\exists\ S'\subseteq S$ such
 that $S_x=S'$ a.e. We assume (without loss of generality) that $S'=S$.\par We claim
 that $\mu(a)>0,\ \forall\
a\in\a_k,\ k\ge 1$. To see this, let $a=[a_1,\dots,a_k]\in\a_k$,
set $S_a:=\{a_1,\ldots,a_k\}$, and let $n_s(a)$ denote the number of times
$s$ appears in $a$. Consider the collection $\Cal K_+^\ast$ of all families of disjoint subsets of $\N$ $\{K_s:\ s\in S_a\}$ such that $|K_s|=n_s(a)$ for all $s$. For every $\Cal K_+\in\Cal K_+^\ast$, set
$$
A(\Cal K_+):=\{x=(x_1,x_2,\dots)\in X:\ x_n=s\ \forall\ n\in K_s,\ s\in
S_a\}.
$$ By recurrence,
$
X\overset{\mod\ \mu}\to{=}[N_s\ge n_s(a)\ \forall\ s\in S_a]=\bigcup_{\Cal K_+\in\Cal K_+^\ast}A(\Cal K_+).
$
Since $\Cal K_+^\ast$ is countable, we can choose $\Cal K_+\in\Cal K_+^\ast$ s.t. $\mu[A(\Cal K_+)]\neq 0$.
Evidently $A(\Cal K_+)\overset{\Cal E}\to\longrightarrow
a$. Since $\mu$ is exchangeable, $\mu(a)>0$.

\par Now let $\b$ be the collection of all cylinders of positive finite
measure.
For each $a\in\b\cap\a_n$, define $\nu_a\in\Cal P(X)$ by
$\nu_a(B):=\tfrac{\mu(a\cap T^{-n}B)}{\mu(a)}$, then $\nu_a$ is
$\Cal E$-invariant, ergodic and by  theorem 2.1 (part 2),
a $(X,T,\a)$-product measure.
\par We claim that $\nu_a$ does not depend on  $a\in\b$.
To see this note that for $a\in\b,\ b\in\a_n$ (some $n\ge 1$) we
have $a b\in\b$ and for $c\in\a_k$,
$$\nu_{a b}(c)=\tfrac{\mu(a b c)}{\mu(a b)}=
 \tfrac{\mu(a)\nu_a(b)\nu_a(c)}{\mu(a)\nu_a(b)}=\nu_a(c).$$
 Thus for $a,b\in\b,\ c\in\a_k$, we have $ab\overset{\Cal
E}\to\longrightarrow ba,\ abc\overset{\Cal
E}\to\longrightarrow bac$, whence
 $$\nu_a(c)=\nu_{a b}(c)=\tfrac{\mu(a b c)}{\mu(a
 b)}=\tfrac{\mu(b a c)}{\mu(b a)}=\nu_{b
 a}(c)=\nu_b(c).$$
Writing $\nu_a=\nu\ \
 (a\in\b)$, we see that for $a,b\in\b$: $$\mu(a)\nu(b)=\mu(ab)=\mu(ba)=\mu(b)\nu(a).$$
   Fixing $\b_0\subset\b$ so that
$X=\biguplus_{b\in\b_0}b\ \mod\mu$ and
 fixing $a\in\b,\ \nu(a)>0$, we have
 $$\mu(X)=\sum_{b\in\b_0}\mu(b)=\sum_{b\in\b_0}\tfrac{\mu(a)\nu(b)}{\nu(a)}\le
 \tfrac{\mu(a)}{\nu(a)}<\infty$$
 and part 1 follows from part 2 of theorem 2.1.

\medskip
In order to prove part 2,  assume by way of contradiction that $\mu$ is locally finite, but not recurrent. Fix $a\in S$ and $N\in\Bbb N$ so that $\mu([N_a=N])>0$. By local finiteness,
$\exists$ a cylinder $b:=[b_1,\dots,b_K]$ with $0<\mu(b\cap [N_a=N])<\infty$ and $\#\{1\le i\le K:\ b_i=a\}=N$.
$\exists\ \s $ a permutation so that  $B:=[b_{\s(1)},\dots,b_{\s(K)}]=[\undersetbrace{N-\text{\rm times}}
\to{a,\dots,a},B_{N+1},\dots,B_K]$. By exchangeability, $0<\mu(B\cap [N_a=N])=\mu(b\cap [N_a=N])<\infty$. By local finiteness, we may
assume  $0<\mu([B_{N+1},\dots,B_K])<\infty$.
 \par For $ n\in\Bbb N$, define $\kappa_n:B\cap [N_a=N]\to [B_{N+1},\dots,B_K]\cap [N_a=N]$ by
 $$\align\kappa_n(\undersetbrace{N-\text{\rm times}}\to
{a,\dots,a},B_{N+1}, &\dots,B_K,x_1,x_2,\dots):=\\ &(B_{N+1},\dots,B_K,x_{1},\dots,x_n,
\undersetbrace{N-\text{\rm times}}\to
{a,\dots,a},x_{n+1},\dots).\endalign$$
But $\mu(\kappa_n(B\cap [N_a=N]))=\mu(B\cap [N_a=N])>0\ \forall\  n\in\Bbb N$, in
contradiction to $0<\mu([B_{N+1},\dots,B_K])<\infty$ as
$[B_{N+1},\dots,B_K]\supseteq\biguplus_{n\geq 1} \kappa_n(B\cap [N_a=N]).$
\hfill\qed

\subheading{Examples}   Consider $X=\{0,1\}^\Bbb N$ with $T$
the shift and $\a=\{[0],[1]\}$.
\roster
\item {\it The mixture above may not
be finite:\/}  Set
$\mu_p:=\prod((1-p)\d_0+p\d_1)$, then the measure
$\mu:=\sum_{n=1}^\infty\tfrac1{n+1}\mu_{\frac1n}$ is a
$\s$-finite, infinite mixture of $(X,T,\a)$-product measures. If
$\nu_n:=n\mu_{\frac1n}$, then
$\mu:=\sum_{n=1}^\infty\tfrac1{n(n+1)}\nu_n$ is an average and
$\mu([1])=\nu_n([1])=1\ \forall\ n\ge 1$.

\item {\it The recurrence condition cannot be removed:\/}
  Let $e_n\in X,\ (e_n)_i=0\ \forall\ i\ne n,\ (e_n)_n=1$,
then $\mu:=\sum_{n=1}^\infty\d_{e_n}\in\goth M(X)$ is $\Cal
E(X,T,\a)$-invariant, ergodic but not a mixture of
$(X,T,\a)$-product measures. This example is topologically
$\s$-finite, but not locally finite, or recurrent.

\item {\it Comparison with the Glimm--Effros theorem:\/}  Let $X_\infty:=\{x\in X:\
N_0(x)=N_1(x)=\infty\}\in\B(X)$. By theorem 2.1, there are
non-atomic  $\Cal E(X,T,\a)\cap (X_\infty\x X_\infty)$-invariant,
ergodic probabilities. Thus, by the Glimm-Effros theorem (above),
there are uncountably many mutually singular, non-atomic,
infinite, $\s$-finite, recurrent, $\Cal E(X,T,\a)$-invariant,
ergodic measures on $X$. Proposition 2.2 says that these measures cannot be
topologically $\s$-finite.
\endroster

\heading Part 2. Exchangeability for topological Markov shifts
\endheading

\heading\S3 Topological Markov shifts and existence of conformal measures
\endheading

\subheading{Topological Markov shifts}
A {\it topological Markov shift} (TMS)
is a fibred system $(X,T,\a)$ which satisfies the {\it Markov property} -- for every $A\in\a$, $T(A)$ is $\a$--measurable -- and whose topology is generated by the set of cylinders.

The elements of $\a$ are called the {\it states} of the shift, and the matrix $(A_{ab})_{\a\x\a}$ with $A_{ab}=1$ if $[a,b]\neq\emptyset$ and $A_{ab}=0$ otherwise is called the {\it transition matrix} of the shift. A TMS with set of states $\a$ and transition matrix $(t_{ab})_{\a\x\a}$ can be canonically identified
with the set
$
\{x=(x_1,x_2,\ldots):\forall i\ t_{x_i x_{i+1}}=1\}
$
together with the action of the left shift and the relative product topology in $\a^\Bbb N$.

 A topological Markov shift is {\it topologically transitive} is for any $a,b\in\a$ there exists some $n$ such that $a\cap T^{-n}b\neq\emptyset$, and {\it topologically mixing} (or simply {\it mixing}) if for all $a,b\in\a$, $a\cap T^{-n}b\neq\emptyset$ for all $n$ sufficiently large. (These definitions agree with the standard definitions of these notions for topological dynamical systems.)

A topological Markov shift
is said to be  {\it almost onto} \cite{ANSS}, if $\forall\ b,c\in\a,\ \exists\
n\ge 1,\ b=a_0,a_1,\dots,a_n=c\in\a$
  such that  $Ta_k\cap Ta_{k+1}\ne\emptyset\ \ (0\le k\le n-1)$. This property comes up in the study of ergodicity done in \S 4.

The following example (example 5.2 in \cite{Pe-S}) shows
that mixing $\nRightarrow$ almost onto. Let
$X\subset\{0,1,2\}^\Bbb N$ be the TMS with transition matrix
$A:=\left(\smallmatrix & 1 & & 1 & & 0 & \\
 & 0 & & 0 & & 1 & \\  & 1 & & 1 & & 0 &
\endsmallmatrix\right)$. Evidently, $(X,T,\a)$ is mixing, but not almost onto.

\subheading{Existence of conformal measures} As we shall see below (theorems 5.0, 5.1, and 5.2) the conformal measures which appear in the study of exchangeability for TMS satisfy the conditions of proposition 2.0, and therefore have $\a$--measurable derivatives. It is therefore enough to clarify what $\a$--measurable functions can appear as derivatives. This is done in the next proposition:

 \proclaim{Proposition 3.0}
Let $(X,T,\a)$ be a mixing TMS with countable state space
 $S$ and let $\pi:S\to\Bbb R^+$.
   There is a $\mu\in{\frak M}_\a(X)$ with
$ \tfrac{d\mu\circ T}{d\mu}\ \a$-measurable and  $\mu([s])=\pi_s\
\ (s\in S)$ $\iff$  $\Pi(s):=\sum_{t\in
S}A_{s,t}\pi_t<\infty\ \forall\ s\in S$. In this case, the measure
$\mu$ is a Markov measure, and is unique.
\endproclaim
\demo{Proof} It is easy to see that  $\tfrac{d\mu\circ T}{d\mu}$ is
$\a$-measurable iff $\mu$ is a Markov measure  of the form
$\mu([s_1,s_2,\dots,s_n])=\pi_{s_1}p_{s_1,s_2}\dots p_{s_{n-1},s_n}$,
where $$p_{s,t}=\frac{A_{s,t}\pi_th(s)}{\pi_s}\text{ and }\sum_{t\in S}A_{s,t}\pi_t=\frac{\pi_s}{h(s)}\
\forall\ s\in S\text{ and
some\ \ }h:S\to\Bbb R_+\ .$$ The first condition implies that
$\frac{d\mu\circ T}{d\mu}(x)=\frac{\pi_{x_2}}{\pi_{x_1}p_{x_1,x_2}}=\frac1{h(x_1)} $; The second guarantees that
$(p_{s,t})$ is a stochastic
matrix.

We prove the equivalence proclaimed above.
\f $(\Rightarrow$) \ \ Suppose that
$\mu([s_1,s_2,\dots,s_n])=\pi_{s_1}p_{s_1,s_2}\dots
p_{s_{n-1},s_n}$ where $p_{s,t}=\frac{A_{s,t}\pi_th(s)}{\pi_s}$
for some $h:S\to\Bbb R_+,$ then
$\Pi(s):=\sum_{t\in S}A_{s,t}\pi_t=\sum_{t\in S}\frac{\pi_sp_{s,t}}{h(s)}=\frac{\pi_s}{h(s)}
<\infty.$ \f $(\Leftarrow$) \ \ Set
$\mu([s_1,s_2,\dots,s_n]):=\pi_{s_1}p_{s_1,s_2}\dots
p_{s_{n-1},s_n}$ where $p_{s,t}=\frac{A_{s,t}\pi_t}{h(s)\pi_s}$
and $h(s):=\tfrac{\Pi(s)}{\pi_s}$. This is a stochastic
matrix: $\sum_{t\in S}p_{s,t}=\sum_{t\in
S}\frac{A_{s,t}\pi_t}{h(s)\pi_s}=\frac{\Pi(s)}{h(s)\pi_s}=1.$
\hfill\qed\enddemo

\noindent
{\em Remarks\/}:
 It is not clear  when the measures in proposition 3.0 are recurrent. However:
\roster
\item We give examples of recurrent measures with $\a$-measurable
derivatives on simple aperiodic random walks below.
\item Let $(X,T,\a)$ be a
mixing TMS with state space $S$. Using the methods of \cite{VJ},
one can characterize those  $h:S\to\Bbb R^+$ for which there is a
recurrent $\nu\in{\frak M}_\a(X)$ with $ \tfrac{d\nu\circ
T}{d\nu}(x)=\tfrac{c}{h(x_0)}$ for some $c>0 $. The condition is
that $\exists\  s\in S$ for which the power series $$F(x):=
\sum_{n=1}^\infty z^n\sum_{x\in [ s],\ T^n x=x}h(x_0)h(x_1)\dots
h(x_{n-1})$$ has a positive radius of convergence $R_s$ and
$F(R_s)=\infty$.
\item As shown in \cite{Sa}, the previous condition always holds
if $\sum_{s\in S}h(s)<\infty$ and $(X,T,\a)$ has the {\it big
images and preimages (BIP)} property: $\exists\
b_1,\dots,b_N\in S$ such that $\forall a\ \exists i,j$ such that
$A_{b_j a}A_{a b_i}=1$.
\endroster

\heading \S 4
$\Cal E(X,T,\a)$--Ergodicity of measures with $\a$--measurable derivatives
\endheading

\subheading{Topological transitivity and aperiodicity}
An equivalence relation on a topological space is called {\em topologically transitive}, if at least one of its equivalence class is dense. This is a necessary condition for  ergodicity w.r.t a globally supported measure. For some classes of measures, it is also a sufficient condition.

We study topological transitivity for $\goth T(T_\phi)\equiv\goth T(T)_{\wh{\phi}}$ in the case of a TMS $(X,T,\a)$ and an $\a$--measurable $\phi:X\to \Bbb G$, where  $\Bbb G$
denotes a locally compact, Abelian, Polish (LCAP) topological
group.
\par Suppose that $\phi:X\to\Bbb G$ is continuous. For $s\in S\
n\ge 1$, let $\Pi_{n,s}:=\{x\in [s]:\ T^nx=x\}$, let
$\Pi_{n}:=\{x\in X:\ T^nx=x\}=\bigcup_{s\in S}\Pi_{n,s}$ and let
$$\Bbb F_{\phi,s}:=\<\{\phi_n(x)-\phi_n(y):\ n\ge
1,\ x,y\in\Pi_{n,s}\}\>,$$
$$\Bbb F_{\phi}:=\<\{\phi_n(x)-\phi_n(y):\ n\ge 1,\
x,y\in\Pi_{n}\}\>,$$ where $\<A\>$ denotes the subgroup
generated by $A$. \par Note that if $\phi$ is $\a$-measurable,
then $\{\phi_n(x)-\phi_n(y):\ n\ge 1,\ x,y\in\Pi_{n,s}\}$ is a
group for every $s$. If in addition $(X,T,\a)$ is mixing, then
$\{\phi_n(x)-\phi_n(y):\ n\ge 1,\ x,y\in\Pi_{n}\}$ is also a
group.
\proclaim{Proposition 4.0}If $(X,T,\a)$ is mixing and
$\phi:X\to\Bbb G$ is $\a$-measurable, then $\Bbb F_{\phi,s}=\Bbb
F_{\phi}\ \forall\ s\in S.$\endproclaim\demo{Proof}  We
first show that $\Bbb F_{\phi,s}=\Bbb F_{\phi,t}\ \forall\ s,t\in
S.$
To see this fix (using mixing) $k\ge 2,\ sat,tbs\in\a_k$. If $g\in\Bbb F_{\phi,s}$, then
 $\exists\ n\ge 1,\ \o{sx},\ \o{sy}\in\Pi_{n,s}$ with $g=\phi_n(\o{sx})-\phi_n(\o{sy})$. It follows that
 $\o{tbsxsa},\ \o{tbsysa}\in\Pi_{t,n+2k-2}$ and
 $$\phi_{n+2k-2}(\o{tbsxsa})-\phi_{n+2k-2}(\o{tbsysa})=\phi_n(\o{sx})-\phi_n(\o{ty})=
 g\in\Bbb F_{\phi,t}.$$

Now we  show that $\Bbb F_{\phi,u}=\Bbb F_{\phi}\ \forall\ u\in S.$
Let $g\in\Bbb F_{\phi}$, then $\exists\, n\ge 1,\ s,t\in S,\ \o{sx}\in\Pi_{n,s},\ \o{ty}\in\Pi_{n,t}$
with $g=\phi_n(\o{sx})-\phi_n(\o{ty})$. By the previous paragraph, it suffices to
show that $g\in\Bbb F_{\phi,s}$. To this end, using mixing,
fix $k\ge 1,\ sat,tbs\in\a_k$, then $\o{sxsatb},\
\o{satytb}\in\Pi_{s,n+2k-2}$ and
$\phi_{n+2k-2}(\o{sxsatb})-\phi_{n+2k-2}(\o{satytb})=\phi_n(\o{sx})-\phi_n(\o{sy})=g\in\Bbb
F_{\phi,s}.$\hfill\qed\enddemo We call the continuous
$\phi:X\to\Bbb G$ {\it (topologically) aperiodic} if $\o{\Bbb
F}_{\phi}=\Bbb G$. It follows from lemma 4.3 (below)  that this is equivalent to the
absence of non-trivial solutions to the functional equation
$\g\circ\phi=\l\tfrac{g\circ T}{g}$ where $\g\in\widehat{\Bbb G},\
\l\in\Bbb S^1$ and $g:X\to\Bbb S^1$ continuous.

\proclaim{Proposition 4.1 (tail
transitivity)}\ \ If $(X,T,\a)$ is mixing and
$\phi:X\to\Bbb G$ is $\a$-measurable and aperiodic, then $\goth
T(T_\phi)$ is topologically transitive on $X\x\Bbb
G$.\endproclaim\demo{Proof}
Fix $g\in\Bbb F_\phi,\ k\ge 1,\ a,b\in\a_k$. We'll show that
$\exists\ a'\subseteq a,\ b'\subseteq b$ so that
$a'\x\{0\}\overset{\goth T(T_\phi)}\to{\rightarrow}b'\x\{g\}$.\par
Indeed, by the mixing of $(X,T,\a), \exists\ \ell\ge 1,\ c,d\in
\a_\ell,\ s\in S$ so that $acs,\ bds\in\a_{k+\ell+1}$. By proposition 4.0,
$\Bbb F_{\phi,s}\owns g$ and $\exists\ p\ge 1,\
x,y\in \Pi_{p,s}$ so that
$\phi_p(y)-\phi_p(x)=g+\phi_{k+\ell}(acs)-\phi_{k+\ell}(bds)$.
\par It follows that
$[acx_1^ps]\x\{0\}\overset{\goth
T(T_\phi)}\to{\rightarrow}[bdy_1^ps]\x\{g\}$. Transitivity now
follows from aperiodicity: $\o{\Bbb F}_{\phi}=\Bbb G$. \hfill\qed\enddemo

By \cite{ADSZ}, $(X,T,\a)$ is almost
onto iff
  $F^\natural:X\to\Bbb Z^{\a\setminus \{a_0\}}_0$ (with $a_0\in\a$  fixed, as defined above)
is aperiodic. It is also shown there that in this case an $\a$-measurable
  $\phi:X\to\Bbb G$ is aperiodic whenever $\<\{\phi(x)-\phi(y):x,y\in X\}\>$ is dense in $\Bbb G$.

 To illustrate this, we give an aperiodicity proof that the example given in \S3 is not almost onto: Define $F^\natural:X\to\Bbb Z^{\{1,2\}}$ by
$F^\natural(x)=e_{x_1}$ if $x_1=1,2$ and $F^\natural(x):=0$ if
$x_1=0$. We see that $F^\natural_1=F^\natural_2\circ T$, whence
$$F^\natural=(F^\natural_2,F^\natural_2)+(F^\natural_2\circ T,0)-(F^\natural_2,0).$$
Therefore $F^\natural:X\to\Bbb Z_0^{\{1,2\}}$ is not aperiodic, and so $X$ cannot be almost onto.

An $\a$--measurable $\phi$ which is not aperiodic, can be modified by a coboundary to be aperiodic as a function into a smaller group. We explain how to do this.

 \proclaim{Livsic cohomology theorem 4.2\ \ \ [L]}
Suppose that $(X,T,\a)$
is a topologically transitive TMS. Let the topology on $\Bbb G$ be
 generated by the  norm
$\|\cdot\|_\Bbb G$ and suppose that $N\ge 1,\ \phi:X\to\Bbb G$ is
$(\a,N)$-H\"older continuous.
If $\phi_n(x)=0\ \forall\ x\in X,\ T^nx=x$, then $\exists\
g:X\to\Bbb G$  $(\a,N)$-H\"older continuous such that
$\phi=g-g\circ T$.\endproclaim
\demo{Proof}  Fix $z\in X$ with
$\o{\{T^nz:\ n\in\Bbb Z\}}=X$, fix $G(z)\in\Bbb G$ and define
$G:\{T^nz:\ n\in\Bbb Z\}\to\Bbb G$ by $G(T^nz):=G(z)+\phi_n(z)$.
We'll extend the domain of definition of $G$ to all $X$ and show
that the extension is $(\a,N)$-H\"older continuous.
\par Suppose that $\rho\in (0,1)$ and
$\|\phi(x)-\phi(y)\|\le M\rho^n\ \forall\ n\ge N,\ x,y\in X,\
x_1^n=y_1^n$. It suffices to show that $$\|G(u)-G(v)\|\le
\tfrac{M\rho^n}{1-\rho}\ \forall\ n\ge N,\ u,v\in \{T^nz:\
n\in\Bbb Z\},\ u_1^n=v_1^n.$$ \par Accordingly, suppose that $n\ge
N,\ \ell
>k,\ (T^kz)_1^n=(T^\ell z)_1^n$ and let $y:=\o{z_{k+1}^{\ell}}$, then
$y\in X,\ T^{\ell-k}y=y$ and $y_1^{\ell-k+n}=z_{k+1}^{\ell+n}$.
Thus
$$\align \|G(T^kz)-G(T^\ell z)\| & =\|\phi_{\ell-k}(T^kz)\|\\ &=
\|\phi_{\ell-k}(T^kz)-\phi_{\ell-k}(y)\|\ \ \ \ \ (\text{ because }
\phi_{\ell-k}(y)=0)\\&\le \sum_{i=0}^{\ell-k-1}\|\phi(T^{k+i}z)-\phi(T^iy)\|\\
&=M\sum_{i=0}^{\ell-k-1}\rho^{n+\ell-k-i}\ \ \ \ \
(\text{ because }z_{k+i+1}^{\ell+n}=y_{i+1}^{\ell-k+n})\\
& \le \frac{M\rho^n}{1-\rho}.\endalign$$
 \hfill\qed\enddemo

\proclaim{Cohomology lemma 4.3\ \ \
c.f.\ [Pa-S]}
If $(X,T,\a)$ is mixing and $\phi:X\to\Bbb G$ is $\a$-measurable,
then
$\phi=a+g-g\circ T+\o\phi$ where
  $a\in\Bbb G,\ \o\phi:X\to \o{\Bbb F}_\phi$  and
$g:X\to\Bbb G$ are both $\a$-measurable such that $\o\phi:X\to
\Bbb F_\phi$ is aperiodic.\endproclaim \demo{Proof}
 By assumption, $\forall\ N\ge 1$ such that $\Pi_N\ne\emptyset,\ \exists\
a_N\in\Bbb G$ so that $\phi_N(x)=a_N\ \mod\Bbb F_\phi\ \forall\
x\in\Pi_N$. Evidently $ka_N=N a_k\ \mod\Bbb F_\phi$ whence
$(\phi_N-a_N)_k=0\ \mod\Bbb F_\phi$ whenever $\Pi_k,\
\Pi_N\ne\emptyset$.\par By Livsic's theorem, $\forall\ N\ge 1$
such that $\Pi_N\ne\emptyset,\ \exists\ g^{(N)}:X\to\Bbb G$
$\a$-measurable such that $$\phi_N=a_N+g^{(N)}-g^{(N)}\circ T\mod
\ \Bbb F_\phi.$$
\par Since $(X,T,\a)$ is mixing, $\exists\ p,q\in\Bbb N$
relatively prime, and $u,v\in X,\ T^pu=u,\ T^qv=v$. Suppose  that
$k,\ell\in\Bbb N$ satisfy $kp-\ell q=1$, then $\mod \ \Bbb
F_\phi$, $$\align  \phi&= \phi_{kp}- \phi_{\ell q}\circ T\\ &=
ka_p-\ell a_q+g^{(p)}-g^{(p)}\circ T^k-g^{(q)}\circ T+g^{(q)}\circ
T^{\ell+1}\\ &= ka_p-\ell a_q+g_k^{(p)}-g_k^{(p)}\circ
T-g_\ell^{(q)}\circ T+g_\ell^{(q)}\circ T^2\\ &=a+g-g\circ
T\endalign$$ where $a:=ka_p-\ell a_q$ and
$g:=g_k^{(p)}-g_\ell^{(q)}\circ T$.
\par Now let $\o\phi:=\phi-(a+g-g\circ T)$, then $\o\phi:X\to\Bbb F_\phi$ is $\a$-measurable and aperiodic since $\<\o\phi_n(x)-\o\phi_n(y):\ n\ge 1,\ x,y\in \Pi_n\>=\Bbb
F_\phi.$ \hfill\qed\enddemo
\subheading{Ergodicity}
Let $(X,T,\a)$ be a TMS and $\mu\in\goth M_\a(X)$ be a
$\goth G(T)$-nonsingular measure. Consider the following properties:
\roster
\item {\it The Gibbs property:\/} $\tfrac{d\mu\circ
T}{d\mu}$ is $(\a,1)$-H\"older continuous;
\item {\it The Markovian Gibbs property:\/} $\tfrac{d\mu\circ T}{d\mu}$ is
$(\a,2)$-H\"older continuous.
\endroster
Some examples: If
$(X,T,\a)$ is full, then any global $(X,T,\a)$-product has the
Gibbs property; and if $(X,T,\a)$ is a TMS , then any global
Markov measure on $X$ has the Markovian Gibbs
property.
\proclaim{Lemma 4.4} Suppose that $(X,T,\a)$ is full
and that  $m\in\Cal P(X)$ is a globally supported measure with the
Gibbs property so that $(X,\B,m,T)$ (where $\B:=\B(X)$) is
conservative and exact, then $m$ is $\Cal
E(X,T,\a)$-ergodic.\endproclaim\demo{Proof}\ \ It is easily
checked  that $(X,\B,T,m,\a)$ is a Gibbs-Markov map in the sense
of \cite{AD}. We have already seen that
$\Cal E(X,T,\a)=\goth T(T,F^\natural)\cong\goth T(T_{F^\natural})\cap (X\x\{0\})^2.$
By construction of $F^\natural$, $\{F^\natural(x)-F^\natural(y):\
x,y\in X\}$ generates $\Bbb Z^{\a\setminus \{a_0\}}_0$ and so
$F^\natural:X\to\Bbb Z^{\a\setminus \{a_0\}}_0$ is aperiodic. It
follows from the theorem in \cite{AD} that $T_{F^\natural}$ is
exact, whence $\goth T(T_{F^\natural})$ is ergodic with respect
to $m\x m_{\Bbb Z^{\a\setminus \{a_0\}}_0}$.
Restricting this
equivalence relation to $(X\x\{0\})^2$, we see that $\Cal E(X,T,\a)$
must also be ergodic with respect to $m$.\hfill\qed\enddemo
\proclaim{
Proposition 4.5} Suppose that $(X,T,\a)$ is a mixing TMS, that  $\Bbb G$ is countable and that
 $\phi:X\to\Bbb G$ is
aperiodic and has finite memory. Let $m\in\goth M_\a(X)$ be a
globally supported  measure with the Markovian Gibbs property such
that $(X,\B,m,T)$ is conservative and exact, then $m\x m_\Bbb G$ is
$\goth T(T)_{\wh{\phi}}=\goth T(T_{\phi})$--ergodic.\endproclaim\demo{Proof} After possibly recoding, we may
assume that $\phi(x)=\phi(x_1)$.  Since $m$ is globally supported and Gibbs,
$m$ is $\goth T(T)$-nonsingular.
\par Let $A=\biguplus_{g\in\Bbb G}A_g\x\{g\} \in\B(X\x\Bbb G)$ be
$\goth T(T_\phi)$-invariant. We claim first that $\forall\ s\in
S,\ g\in\Bbb G$, either $s\subseteq A_g\ \mod m$, or $m(s\cap
A_g)=0$.
\par To see this, fix $s\in\a,\ g\in\Bbb G$ and let
$\widetilde\a_s\subset\B(s)$ be the first return time partition (to $s$), and
$T_s:s\to s$ the induced map. The Markovian Gibbs property
implies that $m|_s$ is a $(s,T_{s},\widetilde\a)$-product measure. Since
$$\Cal E(s,T_{s},\widetilde\a)\subset\Cal E(X,T,\a)\cap (s\x
s)\subset\goth T(T_\phi)\cap (s\x s),$$ we have that $s\cap A_g$
is  $\Cal E(s,T_{s},\widetilde\a)$-invariant. By lemma 4.4, either
$s\subseteq A_g\ \mod m$, or $m(s\cap A_g)=0$. If $m\x m_{\Bbb
G}(A)>0$, then $\exists\ g\in\Bbb G$ with $m(A_g)>0$, whence
$\exists\ s\in S$ with $s\x\{g\}\subset A\ \mod m\x m_{\Bbb G}$.
Since $\phi$ is aperiodic, by proposition 4.1, $\goth T(T_\phi)$
is topologically transitive on $X\x\Bbb G$ and $A=X\x\Bbb G\ \mod
m\x m_{\Bbb G}$. \hfill\qed\enddemo

\noindent
The following corollary was proved for Markov measures in \cite{Gr}.\proclaim{
Corollary 4.6} Suppose that $(X,T,\a)$ is a mixing TMS. Let
$m\in\goth M_\a(X)$ be a globally supported measure with the
Markovian Gibbs property  so that $(X,m,T,\a)$ is conservative and
exact, then
\roster
\item if $(X,T,\a)$ is almost onto, then $m$ is
$\Cal E(X,T,\a)$-ergodic;
\item if not and
$F^\natural=a+g-g\circ T+\o{ F^\natural}$ where $\o{
F^{\natural}}:X\to\Bbb F_{F^\natural}$ is aperiodic, then for each
$b\in\Bbb Z^{S\setminus\{s_0\}}_0,\ m([g\in b+\Bbb
F_{F^\natural}])>0$, $m|_{[g\in b+\Bbb F_{F^\natural}]}$ is
 $\Cal
E(X,T,\a)$-non-singular and ergodic.
\endroster
\endproclaim\demo{Proof}
Fix $s_0\in S$ and apply propositions 4.5, 1.0 and identity (0.2) to the aperiodic
$F^\natural:X\to\Bbb Z^{S\setminus\{s_0\}}_0$ for (1), and to $\o{
F^{\natural}}:X\to\Bbb F_{F^\natural}$ for (2).
\hfill\qed\enddemo
\noindent
Of course this corollary applies to measures with $\a$--measurable derivatives.

\heading{\S5 From  exchangeable measures to conformal measures (TMS)}
\endheading

\subheading{Invariant recurrent measures for cocycle sub-relations} The following results are geared towards showing that recurrent, exchangeable,
ergodic, locally finite measures have $\a$-measurable
derivatives. Noting that $\Cal E(X,T,\a)=\goth T(T,F^\natural)$, we consider the the more general problem of identifying the ergodic invariant measures for $\goth T(T,\phi)$ with $\a$--measurable  $\phi$'s, in the recurrent case. The reader is invited to recall the notation introduced in equation (1.0).

  \proclaim{Theorem 5.0 (globally supported, aperiodic case)}
 Suppose that $(X,T,\a)$ is a mixing TMS, $\Bbb G$ is countable and that
$\phi:X\to\Bbb G$ is $\a$-measurable and aperiodic.
 If $\nu\in\goth M(X )$ is
globally supported, locally finite, $\goth T (T,\phi)
$-invariant, ergodic and
 recurrent,
 then $\nu\in\goth M_\a(X)$ is Markov and
$$
\frac{d\nu\circ T}{d\nu}=ce^{h\circ\phi},\  \text{\rm for some  homomorphism $h:\Bbb G\to\Bbb R$ and $c>0$.}
$$
\endproclaim
\medskip
\noindent
{\it Proof.\/} We divide the proof into several steps.
\medskip
\noindent {\it Step 1.\/} $\nu\in\goth M_\a(X)$, is a $\goth T (T)
$-non-singular Markov measure.

\medskip
To see this fix any $s\in S$ and let $\varphi_s(x):=\min\{n\ge 1:\ T^nx\in [s]\}\ \ (x\in [s])$, $\widetilde\phi:=\sum_{k=0}^{\varphi_s-1}\phi\circ T^k$, and
$\widetilde{\a}=\{[s,a_1,\ldots,a_{n-1},s]: n\in\Bbb N, \forall i\ a_i\in\a\text{ and } a_i\neq s\}$. Note that $\widetilde{\a}$ is a partition of $[s]$ modulo $\nu$ ($\nu$ is recurrent). Also, $\widetilde{\phi}$ is $\widetilde{\a}$--measurable, and it is routine to check that
$$\Cal
E([s],T_{[s]},\widetilde\a_s)\subseteq\goth
G(T_{[s]},(\widetilde\phi,\varphi_s))=\goth
T(T,\phi)\cap([s]\x[s]).$$ By proposition 2.2, part 2,
$\nu|_{[s]}$ is recurrent with respect to
$([s],T_{[s]},\widetilde\a_s)$ and thus, by proposition 2.2, part
1, a mixture of $([s],T_{[s]},\widetilde\a_s)$-products.
\par If $n\ge 1$  and
$[s,a_1,\dots,a_n,t],\ [s,a_1',\dots,a_n',t]$ are cylinders with
$\sum_{j=1}^n\phi(a_i)=\sum_{j=1}^n\phi(a_i')$, then
$$\nu([s,a_1,\dots,a_n,t]\cap T^{-(n+k)}b)=
\nu([s,a_1',\dots,a_n',t]\cap T^{-(n+k)}b)\ \forall\ k\ge 1,\
\text{cylinder}\ b.$$ It follows from this that each component
product of $\nu_{[s]}$ (being a conditional probability with
respect to the tail $\s$-algebra) is $\goth
T(T,\phi)\cap([s]\x[s])$-invariant (c.f. the proof of lemma 12 of
[D-F]). By $\goth T(T,\phi)\cap([s]\x[s])$-ergodicity of
$\nu|_{[s]}$, there is actually only one component, whence
$\nu([s])<\infty$.

To show the Markov property, we note that $\nu|_{[s]}$, being
exchangeable, is partially exchangeable on $[s]$ in the sense of
[D-F]. By proposition 15 of [D-F], it is the restriction to $[s]$
of a Markov measure.  The transition probability of $\nu|_{[s]}$
does not depend on $s\in S$. Being globally supported and Markov,
$\nu$ is  $\goth T (T) $-non-singular.  Step 1 is established.

\medskip
\noindent
{\it Step 2.\/}  Next, let $m\in\goth M(X\x\Bbb G)$ be the unique $\goth T(T_\phi)$-invariant, ergodic measure satisfying
$m(A\x\{0\})=\nu(A)\ \forall\ A\in\B(X)$ (see Proposition 1.0).
We claim that $m\ll\nu\x m_\Bbb G$ and that
$m([s]\x\{g\})<\infty\ \forall\ s\in S,\
g\in\Bbb G$.

\medskip
To prove $m\ll\nu\x m_\Bbb G$, assume
$m(A\x\{g\})>0$. By ergodicity, $\exists\ A'\subset A,\
m(A'\x\{g\})>0$ and $B\in\B(X),\ \nu(B)=m(B\x\{0\})>0$ so that
$A'\x\{g\}\overset{\goth T(T_\phi)}\to\rightarrow B\x\{0\}$, and
in particular $A'\overset{\goth T(T)}\to\rightarrow B$. By $\goth
T (T) $-non-singularity of $\nu,\ \nu(A)\ge\nu(A')>0$.
\par To show that $m([s]\x\{g\})<\infty\ \forall\ s\in S,\
g\in\Bbb G$, note that for each $z\in\Bbb G,\ m\circ Q_z$ is a
$\s$-finite,
$\goth T(T_\phi)$-invariant, ergodic measure where $Q_z(x,y):=(x,y+z)$. Thus either $m\circ Q_z\propto m$, or
$m\circ Q_z\perp m$. Let $$\Bbb H=\Bbb H(m):=\{z\in\Bbb G:\ m\circ Q_z\propto m\},$$ then $\exists$
 a homomorphism $h:\Bbb H\to\Bbb R$ such that $m\circ Q_z=e^{-h(z)}m\ \forall\ z\in\Bbb H$. It follows from this, and $m\ll\nu\x m_\Bbb G$
 that for $A\in\B(X),$
 $$m(A\x\{z\})=\cases & e^{-h(z)}\nu(A)\ \ \ \ z\in\Bbb H,\\ & 0\ \ \ \ z\notin\Bbb H,\endcases$$
 whence the local finiteness of $m$.

\medskip
\noindent
{\it Step 3.\/} $\Bbb H=\Bbb G$.

\medskip
\noindent
Suppose that $g\in\Bbb G$ and fix $s\in S$. By aperiodicity, $\Bbb F_{\phi,s}=\Bbb G$ and
 $\exists\ n\ge 1,\ a=\o{sa_1a_2\dots a_{n-1}},\ b=\o{sb_1b_2\dots b_{n-1}}\in\Pi_{n,s}$ satisfying
 $\phi_n(b)-\phi_n(a)=g$. Since $\nu$ is globally supported,
 $\nu([sa_1a_2\dots a_{n-1}]),\ \nu([sb_1b_2\dots b_{n-1}])>0$.

We claim that
 $[sa_1a_2\dots a_{n-1}]\x\{0\}\overset{\goth T(T_\phi)}\to\rightarrow [sb_1b_2\dots b_{n-1}]\x\{g\}$
 by the map $((s,a_1,a_2,\dots, a_{n-1},x),z)\mapsto ((s,b_1,b_2,\dots, b_{n-1},x),z+g)$. This is because $$(s,a_1,a_2,\dots, a_{n-1},x)\mapsto (s,b_1,b_2,\dots, b_{n-1},x)$$ is a $\goth T(T)$-holonomy
 ($[sa_1a_2\dots a_{n-1}]\to[sb_1b_2\dots b_{n-1}]$) and
 $$\phi_n(s,b_1,b_2,\dots, b_{n-1},x)-\phi_n(s,a_1,a_2,\dots, a_{n-1},x)=\phi_n(b)-\phi_n(a)=g\ \forall\ x\in [s].$$
Thus
 $m(X\x\{g\})>0$ and $g\in\Bbb H$.

\medskip
\noindent
{\it Step 4.\/} We now complete the proof of the theorem by showing that
 $\tfrac{d\nu\circ T}{d\nu}=ce^{h\circ\phi}$ for some $c>0$.

\medskip
The proof of Step 2 shows that $m=e^{-h}\nu\x dm_{\Bbb G}$ (where $m_{\Bbb G}$ is the counting measure on $\Bbb G$). Since $m$ is by definition $\goth T(T_\phi)$--invariant,
for any $\goth T(T)$-holonomy $K$,
 $\tfrac{d\nu\circ K}{d\nu}(x)=e^{h(\widehat\phi(x,Kx))}$.  By proposition 1.1, $\tfrac{d\nu\circ
T}{d\nu}=ce^{h\circ\phi}$ for some $c>0$.
 \hfill\qed

\proclaim{Theorem 5.1\ \ \ (openly supported, periodic case)}
Suppose that $\Bbb G$ is countable, that $(X,T,\a)$ is a mixing TMS, that
$\Bbb F$ is a subgroup of $\Bbb G$  and that
$\phi:X\to\Bbb G$ is $\a$-measurable and of form
$\phi=a+g-g\circ T+\o\phi$ where $a\in\Bbb G,\ g:X\to\Bbb G$ is $\a$-measurable and
$\o\phi:X\to\Bbb F$ is $\a$-measurable and aperiodic.
\par If $\nu\in\goth M(X )$ is locally finite, has clopen support $U$,
is recurrent and $\goth T (T,\phi) $-invariant, ergodic,
 then   $U=[g\in z_0+\Bbb F]$ for some $z_0\in\Bbb G$ and $\exists$ a homomorphism
 $h:\Bbb G\to\Bbb R,\ c>0$ and $\mu\in\goth M_\a(X)$ Markov, such that
$\tfrac{d\mu\circ T}{d\mu}=ce^{h\circ\phi},\ \ \nu=\mu|_{U}.$
\endproclaim
\medskip
\noindent {\it Proof.\/} As in the proof of Step 1 of theorem 5.0,
$\nu\in\goth M_\a(X)$ is the restriction to $U$ of a global Markov
measure and is thus $\goth T(T)\cap (U\x U)$-non-singular.\par
Next, let $m\in\goth M(X\x\Bbb G)$ be the unique $\goth
T(T_\phi)$-invariant, ergodic measure satisfying
$m(A\x\{0\})=\nu(A)\ \forall\ A\in\B(X)$ (see proposition 1.0) and
let $\o m:=m\circ\pi$ where $\pi(x,y):=(x,y-g(x))$. Since
$\pi^{-1}\circ T_\phi\circ\pi=T_{a+\o\phi}$, we have that $\o
m\in\goth M(X\x\Bbb G)$ is $\goth T(T_{\o\phi})$-invariant,
ergodic. Also, $\Bbb H(m)=\Bbb H(\o m)$ (as defined in Step 2 of
the proof of theorem 5.0).

We claim first that $\Bbb H(\o m)\subseteq\Bbb F$.
To see this, assume w.l.o.g. that $m(X\x\{0\})>0$. If $A\in\B(X),\
g\in\Bbb G$ satisfy
 $m(A\x\{g\})>0$, then by ergodicity, $\exists\ A'\subset A,\ m(A'\x\{g\})>0$ and $B\in\B(X),\ m(B\x\{0\})>0$
so that $A'\x\{g\}\overset{\goth T(T_{\o\phi})}\to\rightarrow B\x\{0\}$. Since $\o\phi:X\to\Bbb F$, it follows that
$g\in\Bbb F$.

 Next, we claim that  $\Bbb H(\o m)\supseteq \Bbb F$.
 To see this, fix $s\in S,\ [s]\subset U$ and consider $m_s:=m|_{[s]\x\Bbb G}$ which is
$\goth T(T_\phi)\cap ([s]\x\Bbb G)^2$-invariant, ergodic. Now
$$\goth T(T_\phi)\cap ([s]\x\Bbb G)^2=(\goth T(T)_{\widehat\phi})\cap ([s]\x\Bbb G)^2=
(\goth T(T)\cap([s]\x [s]))_{\widehat\phi}$$ and in a similar
manner to Steps 2 and 3 of the proof of theorem 5.0, we see that
$\Bbb H(\o m)=\Bbb F$.

\par It follows that $\exists\ h:\Bbb F\to\Bbb R$ a homomorphism, $z_0\in\Bbb G$ and a $\s$-finite measure $\o\nu$ on $X$ so that
$$\o m(A\x\{z+z_0\})=\cases & e^{h(z)}\o\nu(A)\ \ \ \ \ \ \ z\in\Bbb
F,\\ & 0\ \ \ \ \ \ \ z\notin\Bbb F.\endcases$$
Assume without loss of generality
 that $z_0=0$. Thus, for $K$ a $\goth T(T)$-holonomy,
$\tfrac{d\o\nu\circ K}{d\o\nu}=e^{-h(\widehat{\o\phi}(x,Kx))}$,
whence as in Step 4 of the of the proof of theorem 5.0,
$\tfrac{d\o\nu\circ
T}{d\o\nu}=ce^{h\circ\o\phi}$ for some $c>0$.
 Now
$$
\align \nu(A) &= m(A\x\{0\})= \o m\circ\pi^{-1}(A\x\{0\})= \o m(\bigcup_{z\in\Bbb G}A\cap
[g=z]\x\{z\})\\
&= \sum_{z\in\Bbb G}\o m(A\cap [g=z]\x\{z\})=\sum_{z\in\Bbb F}e^{h(z)}\o\nu(A\cap [g=z])= \int_{A\cap
[g\in\Bbb F]}e^{h\circ g}d\o\nu.\endalign$$ By theorem 16.1 in
[Fu], $\exists$ a homomorphism $H:\Bbb G\to\Bbb R$ with $H|_\Bbb
F\equiv h$.  Fixing one such $H:\Bbb G\to\Bbb R$, and setting
$d\mu:=e^{H\circ g}d\o\nu$ we see that $\tfrac{d\mu\circ
T}{d\mu}(x)=ce^{H(g(Tx)-g(x))}\tfrac{d\o\nu\circ T}{d\o\nu}
=ce^{H(g(Tx)-g(x))-h(\o\phi(x))}=ce^{-H(\phi(x))}$.\hfill\qed

\subheading{Applications to Exchangeability} Suppose that $(X,T,\a)$ is a  TMS  and that
 $\nu\in\Cal P(X)$ is $\Cal E (X,T,\a) $-invariant and ergodic,
 then
 for $s\in S$ $N_s:=\sum_{n=0}^\infty 1_{[s]}\circ T^n$
is $\Cal E (X,T,\a) $-invariant, whence constant $\nu$-a.e. Recall that $\nu$ is called {\em recurrent}, if $N_s\in\{0,\infty\}$ for all states $s$.
We call $s\in S$ $\nu$-{\it ephemeral} if $1\le N_s<\infty$
$\nu$-a.e., and  the measure $\nu$  {\it ephemeral} if every state
is $\nu$- ephemeral.
Non-ephemerality is strictly weaker than recurrence.

We now apply the previous
results to identify the non-ephemeral ergodic exchangeable measures on topological Markov shifts. For some (though not all) TMS there can be no other topologically $\s$--finite exchangeable measures, see \S6 below.

We maintain the notation of the previous subsection. Let $(X,T,\a)$ be a mixing topological Markov shift, and suppose
that $\mu\in\goth M_\a(X)$ is globally supported and recurrent.

If $\tfrac{d\mu\circ T}{d\mu}$  is $\a$-measurable, then by proposition 2.0, $\mu$ is exchangeable. We claim that its $\Cal E(X)$--ergodic components are given by proposition 4.6. To see that this proposition applies, we need to check that $\mu$ is Markovian, conservative, and exact: The Markov property is clear; The conservativity and exactness of recurrent Markov measures on mixing TMS are well--known (see e.g. [ADU]).

 \par Thus the structure of exchangeable measures with $\a$--measurable derivative is understood. The following two results show that any locally finite, non-ephemeral,
 $\Cal E(X)$- invariant, ergodic measure has a similar form, thus  generalizing and clarifying
 corollary 2.8 in \cite{ANSS1}. We treat the recurrent case separately, because the result is easier to state in this case.

For $\nu\in\goth M(X)$, let $S_\infty=S_\infty(\nu):=\{s\in S:\
 \nu([N_s=\infty])>0\}$, where as before $N_s:=\sum_{n\geq 1}1_s\circ T^n$.
\par Let $X$ be a (shift)-topologically transitive TMS, then  (see e.g. \cite{Ch}),
$X=\biguplus_{k=1}^NX_k$ where $N\in\Bbb N$ and $X_1,\dots,X_N$
are disjoint, clopen subsets of  $X$ with $TX_k=X_{k+1\ \mod N}$;
and each $(X_k,T^N,\a_N)$ is mixing. This decomposition is called
the {\it periodic decomposition} of $X$, $N=N_X$ is called the
{\it period} of $X$ and each $X_k$ is called a {\it basic, mixing}
set for $X$. \proclaim{Theorem 5.2\ \ \ (recurrent case)} Suppose
$(X,T,\a)$ is a  TMS and let $\nu\in\goth M(X )$ is locally
finite, recurrent and $\Cal E (X,T,\a) $-invariant, ergodic, then
there are \roster\item $X'\subset X\cap S_\infty^{\Bbb N}$,\ \
$\Cal E (X,T,\a) $-invariant, such that $(X',T,\a)$ is a
topologically transitive TMS and ;\item a clopen, $\Cal E (X,T,\a)
$-invariant subset $U$ of a basic mixing set for $X'$;
\endroster
so that $\nu=\mu|_{U}$  where $\mu\in\goth M_\a( X')$ is Markov
with $\tfrac{d\mu\circ T^{N_{X'}}}{d\mu}$ $\a$-measurable.
\endproclaim
\demo{Proof} As in the proof of Step 1 of theorem 5.0,  $\nu$ is
the restriction of a Markov measure to a union of initial states.

{ \par The associated stochastic matrix $p$ is recurrent (by
assumption) and irreducible (as $\nu$ is $\Cal E (X,T,\a)
$-ergodic). Thus $X':=\{x\in S_\infty^\Bbb N:\ p_{x_n,x_{n+1}}>0\
\forall\ n\ge 1\}$   is topologically transitive (with respect to
the shift).
\par In case $(X',T,\a)$ is mixing, the result follows from
theorem 5.1.
\par In general $(X',T,\a)$ is transitive with periodic decomposition $X'=\biguplus_{k=1}^NX_k'$.
 Each
$(X_k',T^N,\a_N)$ is mixing and $\goth T(T)=\goth
T(T^N)$-invariant whence $\nu$ is supported on some
$X_{k_0}'$.\par   Thus, $\nu$ is $\Cal E (X_{k_0}',T^N,\a_N)
$-invariant, whence (by the mixing case) a mixture of conformal,
Markov measures. However $\nu$ is $\goth T(T^N)$-ergodic and so
this mixture is trivial. }\hfill\qed\enddemo

\noindent
{\em Remark\/}:\ \ {  The TMS $X'$ is not necessarily mixing, even if $X$ is mixing.
 To see this, let $X\subset\{0,1\}^\Bbb N$ be the (mixing) Fibonacci shift with transition matrix
 $\left(\smallmatrix & 1 & & 1 &  \\
  & 1 & & 0 & \endsmallmatrix\right)$. The measure $\nu:=\d_{(1010\dots)}$ is exchangeable,
  the corresponding TMS
  $X'=\{(1010\dots),\ (0101\dots)\}$ being non-mixing (having period 2)  with transition matrix
 $\left(\smallmatrix & 0 & & 1 &  \\
  & 1 & & 0 & \endsmallmatrix\right)$.}

\proclaim{Theorem 5.3\ \ \ (non-ephemeral case)} If $(X,T,\a)$ is
a  TMS and $\nu\in\goth M(X )$ is locally finite, $\Cal E (X,T,\a)
$-invariant, ergodic and non-ephemeral, then there exist
\roster\item a cylinder $f=[f_1,\dots,f_K]\subset\Sigma$ with
$f_1,\dots f_K$ $\nu$- ephemeral;\item $X'\subset X\cap
S_\infty^{\Bbb N}$ a $\Cal E(X,T,\a)$-invariant  TMS;\item a
clopen, $\Cal E (X,T,\a) $-invariant subset $U$ of a basic mixing
set for $X'$;
\item a  Markov measure $\mu\in\goth M_\a(X')$
\endroster
so that $\tfrac{d\mu\circ T^{N_{X'}}}{d\mu}$ is $\a$-measurable
and
$$\nu=c\sum_{\s\in S_K,\ \s f\cap U\ne\emptyset} \d_{\s
f}\x\mu|_{T^K(\s f)\cap U}$$  where $\s
f:=[f_{\s(1)},\dots,f_{\s(K)}]$.
\endproclaim

\medskip
\noindent
{\it Proof.\/}
 Let $S_e$ be the set of $\nu$--ephemeral states and  $S_\infty$ be the set
of $\nu$--recurrent  states. Let $\goth n(x):=\min\,\{n\ge 1:\ x_n\in S_\infty\}$.
We claim that $\goth n$ is constant and that $\goth n=\sum_{f\in
 S_e}N_f<\infty$. \par To see this, we claim first that for $\nu$-a.e. $x,\ x_k\in S_\infty\ \forall\ k\ge\ \goth n(x)$.
\par Since $S_\infty\ne\emptyset,\ \goth n<\infty\ \nu$-a.e. Suppose that
$x\in X,\ \goth n(x)=K,\ x_K=s\in S_\infty$ and
$\nu([x_1,\dots,x_K])>0$. It follows that $\nu([x_1,\dots,x_K]\cap
T^{-K}(\cdot))$ is a locally finite, $\Cal E
([s],T_{[s]},\widetilde\a)$-invariant, $\goth
T([s],T_{[s]},\widetilde\a)$-ergodic measure on $[s]$ which is, by
proposition 2.2,  a multiple of a
$([s],T_{[s]},\widetilde\a)$-product measure. Thus
$$x_n\in \cases & S_e\ \ \ 1\le n\le\goth n(x),\\ &S_\infty\ \ \ \ \  n\ >\goth n(x),\endcases$$
whence  $\goth n(x)=\sum_{f\in
 S_e}N_f(x)<\infty$, which is constant by $\Cal E (X,T,\a) $-ergodicity.
 The result follows from this, and theorem 5.2.
 \hfill\qed

\heading{\S6 Examples of ephemeral and non-ephemeral exchangeable measures}\endheading

\subheading{Non-ephemeral exchangeable measures with ephemeral states}
The following example is taken from \cite{Pe-S}. Consider the TMS $X$ on the states
$\{0,1,2,3\}^{\Bbb N}$ with transitions
$$\underset{\circlearrowright}\to{0}\leftrightarrow\underset{\circlearrowright}\to{1}{\rightarrow 2\rightarrow 3}\rightarrow 1$$
and the exchangeable $\nu\in\Cal P(X)$ given by
$$\nu([3,\e_1,\e_2,\dots,\e_n]):=\cases & \tfrac1{2^n}\ \ \ \e_1,\e_2,\dots,\e_n=0,1,\\ & 0\ \ \ \ \text{ else.}\endcases$$

\subheading{TMS without locally finite  ephemeral exchangeable measures }
A TMS $(X,T,\a)$ with
state space $S$ and transition matrix $A=(t_{ij})_{i,j\in S\x S}$ has
the {\it finite images property}, if
$A$ has finitely many rows (equivalently $\{Ta:a\in\a\}$ is
finite).
Examples include any  TMS  with  finite state space, the full shift on a countable alphabet, and any shift
obtained from the full shift by removing a finite collection of
edges.

\proclaim{Proposition 6.0} If $(X,T,\a)$ has the
finite images property and $\nu\in\goth M(X)$ is topologically
$\s$-finite, $\Cal E(X,T,\a)$-invariant, ergodic, then $\nu$ is
not ephemeral.
\endproclaim
\demo{Proof}  Let $\nu\in\goth M(X)$ is  topologically
$\s$-finite, $\Cal E(X,T,\a)$-invariant, ergodic. Suppose also
that $\nu\in\goth M(X)$ is ephemeral.\par  For $a\in S$ define
$R(a), C(a):S\to\{0,1\}$ by $ R(a)(b):=t_{ab},\ C(a)(b):=t_{ba}$.
The finite images property is that $\Cal R:=\{R(a):a\in S\}$ is
finite, say equal to $\{R_1,\ldots,R_N\}$.
We claim that $\Cal C:=\{C(a):a\in S\}$ is finite as well. To see this, note that
$$
S=\biguplus_{i=1}^N A_i,\text{\rm where } A_i=\left\{ a\in
S:R(a)=R_i\right\}.
$$
For each $i$, $C(a)$ is constant on $A_i$, because
$$ C(a)(b)=t_{ba}=R(b)(a)=R_i(a)\ \text{\rm  is independent of }\ b\in A_i.$$ It
follows that $|\{C(a):a\in S\}|\leq 2^N$.

\par Fix  $(C,R)\in \Cal C\x\Cal R$ and let $S_{(C,R)}:=\{a\in S:\ C(a)=C,\ R(a)=R\}$.
If $a\in S_{(C,R)},\ x_1,x_2,\dots,x_{n-1}\notin S_{(C,R)}$ and
$[S_{(C,R)}]\cap T^n[a,x_1,\dots,x_{n-1}]\ne \emptyset$ (where
$[S_{(C,R)}]:=\bigcup_{s\in S_{(C,R)}}[s]$), then $\exists\ b\in
S_{(C,R)}$ with $t_{x_{n-1},b}=1$ whence $C(x_{n-1})=1$ and $
T^n[a,x_1,\dots,x_{n-1}]\supset [S_{(C,R)}]$.
 \par It follows that $(A_{(C,R)},T_{A_{(C,R)}},\a_{(C,R)})$ is a full fibred system where
$A_{(C,R)}:=\{x\in [S_{(C,R)}]:\  x_n\in S_{(C,R)}\ \text{\rm
infinitely often}\}$ and
$$\a_{(C,R)}:=\{[a,x_1,\dots,x_{n-1}]:\ n\ge 1,\ a\in S_{(C,R)},\
x_k\notin S_{(C,R)},\ \sum_{b\in S_{(C,R)}}t_{x_{n-1},b}>0 \}.$$
Also
$\Cal E(A_{(C,R)},T_{A_{(C,R)}},\a_{(C,R)})\subset \Cal E(X,T,\a)\cap A_{(C,R)}\x A_{(C,R)}.$
\par  Since
$\Cal C\x\Cal R$ is finite, $\exists\ N_0\ge 1,\ \ (C,R)\in\Cal
C\x\Cal R$ such that $\nu(T^{-N_0}A_{(C,R)})>0$.
By topological $\s$- finiteness, $\exists N\ge N_0,\ a\in\a_N$ such that
$0<\nu(a\cap T^{-N}A_{(C,R)})<\infty$.\par The probability
$q\in\Cal P(A_{(C,R)})$ defined by $q(B):=\tfrac{\nu(a\cap
T^{-N}B)}{\nu(a\cap T^{-N}A_{(C,R)})}$ is  \f$\Cal
E(A_{(C,R)},T_{A_{(C,R)}},\a_{(C,R)})$-invariant and $\Cal
E(X,T,\a)\cap A_{(C,R)}\x A_{(C,R)}$-ergodic.  By theorem 2.1,
part 2, $q$ is a mixture of
$(A_{(C,R)},T_{A_{(C,R)}},\a_{(C,R)})$-product measures, but by
$\Cal E(X,T,\a)\cap A_{(C,R)}\x A_{(C,R)}$-ergodicity, there is
only one component and $q$ is a
$(A_{(C,R)},T_{A_{(C,R)}},\a_{(C,R)})$-product measure.
\par In particular $q$ is recurrent with respect to $(A_{(C,R)},T_{A_{(C,R)}},\a_{(C,R)})$
and $\nu$ cannot be ephemeral. \hfill\qed\enddemo
\subheading{TMS with recurrent and ephemeral locally finite exchangeable measures}
Consider the {\it simple, aperiodic random walks} TMS's $X_d$ for $d\in\Bbb N$ where:
$$X_d:=\{x\in(\Bbb Z^d)^{\Bbb N}:\ \|x_{n+1}-x_n\|_\infty\le 1\ \forall\ n\ge
1\}$$ equipped with the shift map $T$ and partition $\a:=\{[s]:\
s\in\Bbb Z^d\}$. Note that
\roster
\item $(X_d,T,\a)$ is almost onto, so
by proposition 4.6, any  globally supported, Markov measure on $X$ is
$\Cal E(X,T,\a)$-ergodic;
\item $X_d=X^d$ where $X:=X_1$.
\endroster
Now let $\mu\in\goth M_\a(X)$ be  a recurrent,  $\Cal
E(X,T,\a)$-invariant, ergodic, measure on $X$. By theorem 4.3,
$\mu$ is  Markov with $\a$-measurable derivative $\tfrac{d\mu\circ
T}{d\mu}$.
\par The symmetric, global, random walk measure $m$ on
$X$ given by
$$m([a_1,a_2,\dots,a_n])=\cases & \tfrac1{3^n}\ \ \ \ \ \ \
|a_{k+1}-a_k|\le 1,\ 1\le k\le n-1,\\ & 0\ \ \ \ \ \ \ \text{\rm
else}\endcases$$ is null recurrent and of this form. Note however
that $m_d:=\undersetbrace{\text{\rm $d$-times}}\to{m\x\dots\x m}$
is transient on $X_d\ \forall\ d\ge 3$.
\proclaim{Proposition
6.1} For each $d\ge 1$ there is a positively recurrent,
 exchangeable, Markov  measure $\mu\in\goth M_\a(X_d)$.
\endproclaim\demo{Proof} It suffices to find a positively recurrent,
 exchangeable, Markov  measure $\mu\in\goth M_\a(X)$ for then
 $\mu\undersetbrace{\text{\rm $d$-times}}\to{\x\dots\x}\mu\in\goth M_\a(X_d)$
 is as required.\par To this end, fix $0<z<1$, set $\pi_s=z^{|s|}$ and let $\mu\in\goth M_\a(X)$ be
 the Markov  measure with $\mu([s])=\pi_s$ and $\tfrac{d\mu\circ T}{d\mu}$ $\a$-measurable
 (as in proposition 3.0). The underlying stochastic matrix is given by
 for $s\in\Bbb Z,\ i=0,\pm 1$:
$$p_{s,s+i}=\cases & \tfrac{z^{i}}{z^{-1}+1+z}\ \ \ \ s\in\Bbb Z,\ s>0,\ i=0,\pm 1\\ &
\tfrac{z^{-i}}{z^{-1}+1+z}\ \ \ \ s\in\Bbb Z,\ s<0,\ i=0,\pm 1,\\
& \tfrac{z^{|i|}}{2z+1}\ \ \ s=0,\ i=0,\pm 1.\endcases$$
An
invariant distribution on $\Bbb Z$ for $P$ is given by
$$c_t\ =\ \cases\ &\ \tfrac1{1+z+z^2}\ \ \ \ \ t=0,\\ &\
\tfrac{z^{2|t|-1}}{2z+1}\ \ \ \ \ \ t\ne 0.\endcases$$ Since
$\sum_{t\in S}c_t<\infty$, the stochastic matrix $P$ is positively
recurrent.\hfill\qed\enddemo

\proclaim{Proposition 6.2} Set $Z:=\{x\in X:N_s(x)=0\ \forall\/
s\le 0,\ N_s(x)=2\ \forall\/ s\ge 1\}$. There is a unique $\Cal
E(X,T,\a)$-invariant, ergodic measure $\o\nu\in\goth M_\a(X)$
which is carried by $Z$ and  such that $\o\nu([1])=1$. This
measure is non-atomic and  ephemeral.
\endproclaim
\medskip
\noindent
{\it Proof.\/} Set $\Cal E:=\Cal E(X,T,\a)$, $\Om=\{0,1\}^\Bbb N$, $S:\Om\to\Om$ the left shift, and
define $\vartheta:\Om\to\{1,2\}^\Bbb N$ by $\vartheta(x):=x_1$.

\medskip
\noindent
{\it Step 1.\/} Construction of a homeomorphism $\Phi:Z\cap [1]\to \Om:=\{1,2\}^\Bbb N$ which carries $\Cal E\cap (Z\cap [1])^2$ onto $\goth G(S,\vartheta)\cap (\Om\x\Om)$.

\medskip
Call $K\in\Bbb N$
{\it good for $x\in Z\cap [1]$} if $\#\{1\le k\le 2K:\ x_k=s\}=2\ \forall\ 1\le
s\le K$ and let $\kappa(x):=\min\{K\ge 1:\ K\ \text{\rm good for}\
\ x\}$. The possibilities for $\kappa(x)$ are:
\roster
\item
$\kappa(x)=1$ in case $x=(1,1,\tau(x)+(1,1,1,\dots))$ where
$\tau(x)\in Z\cap [1]$, and
\item $\kappa(x)=2$ in case
$x=(1,2,1,2,\tau(x)+(2,2,2,\dots))$ where $\tau(x)\in Z\cap [1]$.
\endroster
Note that $\tau:Z\cap [1]\to Z\cap [1]$  and $\kappa$ are related  by
$$\tau(x):=T^{2\kappa(x)}(x)-(\kappa(x),\kappa(x),\dots).$$
It follows that
$$\align & (x,y)\in\Cal E(X,T,\a)\cap (Z\cap [1]\x Z\cap [1])\ \Leftrightarrow\\ &
 \exists\ k,\ell\ge 1\ \ni\
\kappa(\tau^{j+k}x)=\kappa(\tau^{j+\ell}y)\ \forall j\ge 1, \&\
\sum_{j=0}^k\kappa(\tau^jx)=\sum_{j=0}^\ell\kappa(\tau^jy).\endalign$$
Defining $\Phi:Z\cap [1]\to \Om:=\{1,2\}^{\Bbb Z_+}$ by
$\Phi(x)_n:=\kappa(\tau^nx)$, we have that
$$\Phi\x\Phi\big(\Cal E(X,T,\a)\cap (Z\cap [1]\x
Z\cap [1])\big)= \goth G(S,\vartheta)\cap (\Phi(Z\cap
[1])\x\Phi(Z\cap [1]))$$ where $S$ is the shift map on $\Om$ and
$\vartheta(x):=x_1$.
\par Evidently $\Phi:Z\cap [1]\to \Phi(Z\cap [1])$ is a homeomorphism.
We claim next that $\Phi(Z\cap [1])=\Om$.
To see this set $a(1):=(1,1),\
a(2):=(1,2,1,2)$ and define $\pi:\Om\to\ \Bbb Z^\Bbb N$ by
$$\pi(\om):=(a(\om_0)+s_0,a(\om_1)+s_1,\dots,a(\om_n)+s_n,\dots)$$
where $s_n=s_n(\om)$ is defined by $s_0=0,\ s_{n+1}=s_n+\om_n$.
Evidently $\pi\circ\Phi=\text{\rm Id.}|_{Z\cap [1]}$, whence
$\pi(\Om)=Z\cap [1]$.

\medskip
\noindent
{\it Step 2.\/} There is a unique $\goth G(S,\vartheta)$-invariant $\mu\in\Cal
P(\Om)$.

\medskip
It is not hard to check that if $p(\om)=p^\om\ \
(\om=1,2)$ where $p+p^2=1\  (p=\tfrac{\sqrt 5-1}2)$, and
$\mu\in\Cal P(\Om),\
\mu([\om_1,\dots,\om_n])=\prod_{k=1}^np(\om_k)=p^{\vartheta_n(\om)}$,
then $\mu$ is $\goth G(S,\vartheta)$-invariant.
To see uniqueness of this $\mu$, we note first that
$$\Cal E(\Om)\ \subset\ \goth G(S,\vartheta)\ \subset\ \goth
G(S).$$ If $\o\mu\in\Cal P(\Om)$ is $\goth
G(S,\vartheta)$-invariant, ergodic, then it is $\Cal
E(\Om)$-invariant, whence by de-Finetti an average of product
measures; and $\goth G(S)$-ergodic, whence a product measure.
Writing $\o\mu=\prod q$, we have by $\goth
G(S,\vartheta)$-invariance that $q(1)^2=\mu([11])=\mu([2])=q(2)$
whence $\o\mu=\mu$.

\medskip
\noindent
{\it Step 3.\/} There is
a unique $\Cal E\cap (Z\cap [1])^2$
ergodic probability measure $\nu\in\Cal P(Z\cap [1])$ and this
measure is non-atomic.

\medskip
This follows from Steps 1 and 2: Any $P\in\Cal P(\Om)$
is $\goth T(\Om,S)$-invariant, ergodic iff $P\circ\Phi\in\Cal
P(Z\cap [1])$ is $\Cal E\cap (Z\cap [1])^2$-invariant, ergodic. The required probability is
$\nu:=\mu\circ\Phi$ where $\mu$ is as in Step 2.

\medskip
\par To complete the proof of the theorem, we take for $\o\nu$
the unique $\Cal E$-invariant, ergodic measure $\o\nu$ on
$Z$ so that $\o\nu|_{Z\cap [1]}=\nu$. If $\lambda$ is another
$\Cal E$-invariant, ergodic measure on $Z$ so that
$\lambda(Z\cap [1])=1$, then by Step 3, $\lambda|_{Z\cap [1]}=\nu$,
whence $\lambda=\o\nu$.

Finally, we note that
$Z=\bigcup_{k=1}^\infty
\{(k,k-1,\dots,2,1,1,2,\dots,k-1,k,u^k(x)):\ x\in Z\cap [1]\}$
where $u(x)_n:=x_n+1$, whence $\o\nu\in\goth
M_\a(X)$. In particular, $\o\nu$ is locally finite.\hfill\qed

\noindent
{\em Remark\/}: It follows from
theorem 2.2 in [ANSS] that the only partition-bounded, $\goth
T(T)$-invariant, ergodic measures on $X$ are the  random walk
measures $m_f$ of form
$$m_f([s_1,s_2,\dots,s_n])=\pi_{s_1}p_{s_1,s_2}\dots
p_{s_{n-1},s_n}$$
   with $p_{s,s+\e}:=f_\e\ \ \ (\e=\pm 1)$,
where $f_\e=z^\e p$ for some $0<p\le\tfrac13,\
z=z_{\pm}(p)=\tfrac12(\tfrac1p-1\pm\sqrt{(\tfrac1p-1)^2-4})$ and
$\pi_s=z^s\ \ (s\in\Bbb Z)$ . The only recurrent measure of this
form is the symmetric one with $p=\tfrac13,\ z=1$. All the others
are ephemeral and not $\Cal E(X,T,\a)$- ergodic. Their ergodic
decompositions may be of interest.

\heading
Part 3. Exchangeability for $\b$--expansions
\endheading

\heading{\S7 $\beta$-expansions}\endheading
\subheading{Definition of $\b$--expansions}
Fix a non-integer $\b>1$  and let $I:=[0,1]$. The $\b$--transformation is $T=T_\b:I\to I$, where $T x:=\{\b
x\}$. The {\it $\b$-expansion} of $x\in I$ is
$\pi(x)=\pi_\b(x)=(\e_1,\e_2,\dots)\in \{0,1,\dots,[\b]\}^\Bbb N$
where $\e_n:=[\b T_\b ^{n-1}x]$. Evidently
$$
x=\sum_{n=1}^\infty\frac{\e_n}{\b^n}.
$$
This should be thought of as an expansion to a `non-integer base' (of course, if $\b\in\N\setminus\{1\}$, this is just the expansion to the base $\b$). The $\b$--expansions were introduced by Renyi [Re], who together with Parry [Pa] also initiated the study of their stochastic behaviour (when $x\in [0,1]$ is distributed according to the absolutely continuous invariant probability  measure of $T_\b$).

For $\b$ an integer the collection of $\b$--expansions is a full-shift on $\b=[\b]$ symbols,  and $\e_n$ are i.i.d.'s. For general $\b$'s the collection of $\b$--expansions is not a topological Markov shift (see below), so the digit process $\{\e_n\}$ is not even Markov.

Nevertheless,  various authors have shown that the digit process shares many of the properties of an i.i.d. process: SLLN, Kolmogorov and Hewitt-Savage zero-one laws, CLT, LLT (see [Re], [Pa], [ADSZ]). We shall contribute to this list, by establishing a suitable version of the de-Finetti theorem.

\subheading{Basic properties of the $\b$--transformation} For the  basic results
 see [Ge], [Pa] and [Re]. Here we mention only a few that are used in the sequel.
The $\b$--transformation
 leaves Lebesgue measure $m$ quasi invariant ($m\circ
T^{-1}\sim m$).  It is known that $(I,\B,m,T)$ is an exact endomorphism
with a Lebesgue-equivalent invariant probability.
The triple $(I,T,\a)$ where
$\a:=\{[\tfrac{j}\b,\tfrac{j+1}\b)\}_{j=0}^{[\b]-1}\cup\{[\tfrac{[\b]}\b,1)\}$ is a fibred system.

The {\it $\b$-shift} is the closure of the collection of $\b$--expansions: $X_\b:=\o{\pi_\b(I)}$.
Set
$$\om=\om_\b:=\cases & (\o{\eta_1,\eta_2,\dots,\eta_{q-1},\eta_q-1})\ \ \ \
\pi_\b(1)=(\eta_1,\eta_2,\dots,\eta_{q-1},\eta_q,\o 0),\\ &
\pi_\b(1)\ \ \ \ \text{else.}\endcases$$
The following is in \cite{Pa}:
$$X_\b=\{y\in\{0,1,\dots,[\b]\}^\Bbb N:\ y_k^\infty\prec \om\
\forall\ k\ge 1\}\tag{7.1}$$ where  $x\prec y$ means $\exists\
n\ge 1,\ x_n<y_n,\ \ x_1^{n-1}=y_1^{n-1}$.

\subheading{Full cylinders} A cylinder $[x_1,\dots,x_N]$ with
$T^N[x_1,\dots,x_N]=X_\b$ is called {\it full}. Full cylinders
were considered in [Re], [Pa], [Sm], [I-Ta], [W] and [Bl].

Not every cylinder is full. The image of a general cylinder $a=[a_1,\dots,a_N]\in\a_N$ is given by
$$T^N[a_1,\dots,a_N]=\{y\in X_\b:\ y\prec
\om_{K_N(a)+1}^\infty\},\tag{7.2}$$
where  $$K_N(a):=\cases 0 & \nexists\ 1\le n\le N,\ a_{N-n+1}^N=\om_1^{n},\\
\max\{1\le n\le N,\ a_{N-n+1}^N=\om_1^{n}\}
&
\text{\rm else}.\endcases$$
By (7.1), if $K_N(a)=0$, then $a$ is full.

It is standard to check, noting that  $\om\in X_\b$ and using (7.1), that
if $\nexists n$ s.t. $x_n^\infty=\om$, then $K_N(x_1,\ldots,x_N)=0$ infinitely often.

Set $\G:=\{x\in X_\b:\ \exists\ n\ge 1,\ x_n^\infty=\om\}$, then
$\pi:I\setminus\pi^{-1}\G\to X_\b\setminus\G$ is invertible and
$\pi\circ T_\b=T\circ\pi$ where $T:X_\b\to X_\b$ is the shift.
Accordingly set
$$
X_{\b,0}:=X_\b\setminus\G.
$$
The previous paragraph says that any $x\in X_{\b,0}$ belongs to infinitely many full cylinders.

This allows us to make the following definitions:
$$
\psi:X_\b\setminus\G\to \Bbb N,\  \psi(x):=\min\,\{N\ge 1,\ T^N[x_1,\dots,x_N]=X_\b\},$$
$$
S:X_{\b,0}\to X_{\b,0},\
S(x):=T^{\psi(x)}(x).$$
We call $S$ a {\it
Bernoulli jump transformation} (see page 133 of [Schw]).
Evidently $(X_{\b,0},S,\widetilde\a)$ is a
full fibred system where $$\widetilde\a:=\{a\in\a_n:\ n\ge 1, \
a\subset [\psi=n]\}.$$ We have $$\Cal
E(X_{\b,0},S,\widetilde\a)\subset\goth T(T|_{X_{\b,0}})
\subset\goth T(S)\subset
\goth
G(T|_{X_{\b,0}}).\tag{7.3}$$

\subheading{Ergodic properties of the tail and exchangeable relations of $X_\b$}We abuse notation and denote the Lebesgue measure on $I$ and the measure it induces on $X_\b$ by the same symbol $m$.

\proclaim{Proposition 7.0} For every $\b>1$, $m$ is $\Cal E(I,T,\a)$--ergodic and invariant.
\endproclaim\demo{Proof} Viewing $T$ on $I$, we see that $\frac{dm\circ T}{dm}=T'$ is $\a$--measurable, so $m$ is $\Cal E (I,T,\a)$--invariant by proposition 2.0.

Viewing  $S$ on $I$, we see that it is a piecewise onto affine map,
and so $m$ is a $(S,\widetilde\a)$-product measure. Thus, by
theorem 2.1, $\goth{I}(\Cal
E(X_\b,S,\widetilde\a))\overset{m}\to{=}\{\emptyset,X_\b\}$. By
(7.3) $\goth{I}(\Cal E(X_\b,T,\a))\overset{}\to{\subseteq}\goth{I}(\Cal
E(X_\b,S,\widetilde\a))\overset{m}\to{=}\{\emptyset,X_\b\}$.\hfill\qed\enddemo

\proclaim{Proposition 7.1} $\goth
T(T)$ is uniquely ergodic: If $\mu\in\goth M(X_\b)$ is topologically $\s$-finite and $\goth
T(T)$-invariant, then $\mu=cm$ for some $c>0$.\endproclaim
\demo{Proof} Let  $\mu\in\goth M(X_\b)$ be topologically $\s$- finite and  $\goth
T(T)$-ergodic, invariant. It is easy to check that $\goth
T(T)$-equivalence classes are dense, therefore by local finiteness $\mu$ must be non-atomic. In particular  $\mu(X_\b\setminus X_{\b,0})=0$, so we may work on $X_{\b,0}$.
\par We claim that $\mu$ is locally finite on $X_{\b,0}$. To see this, note first that by topological $\s$-finiteness, there is a cylinder set
of positive, finite measure. By non-atomicity,  we may assume that the cylinder is full.
Namely: $\exists\ b\in\a_N\cap\widetilde\a_{N'}$ with $0<\mu(b)<\infty$. Now let $x\in X$. Since $T^Nb=X$, we may define  the $\goth
T(T)$-holonomy $\kappa:x_1^N\to b$ by $\kappa(x_1^N,z):=(b,z)$. By $\goth
T(T)$-invariance,
$$\mu([x_1^N])\le \mu(\kappa [x_1^N])\le\mu(b)<\infty$$
and $\mu$ is locally finite on $X_{\b,0}$.

By (7.3), $\mu$ is $\goth
T(S)$-ergodic and  $\Cal
E(X_{\b,0},S,\widetilde\a)$-invariant.  By proposition 2.2, $\mu$
is a $(S,\widetilde\a)$-product measure. \par The $\goth
T(T)$-invariance implies that $\exists\ t>0$ so that
$\mu(a)=t^n\ \forall\ a\in\widetilde{\a}\cap\a_n$. To see this,
suppose that $a\in\widetilde\a\cap\a_k,\
b\in\widetilde\a\cap\a_\ell$, then $$[\undersetbrace{\text{\rm
$\ell$-times}}\to{a,\dots,a}]\overset{\goth
T(T)}\to\longrightarrow\ [\undersetbrace{\text{\rm
$k$-times}}\to{b,\dots,b}]\ \ \ \ \text{\rm by }\ \
(\undersetbrace{\text{\rm $\ell$-times}}\to{a,\dots,a},x)\
\mapsto\ (\undersetbrace{\text{\rm $k$-times}}\to{b,\dots,b},x)$$
so $\mu(a)^\ell=\mu(b)^k$. The case $k=\ell$ shows that $\exists\
t_k>0$ so that $\mu(a)=t_k\ \forall\ a\in\widetilde\a\cap\a_k$ and
the other cases show that $\exists\ t>0$ so that $t_k^{\frac1k}=t\
\ \forall\ k\ge 1,\ \widetilde\a\cap\a_k\ne\emptyset$.
This $t>0$ is uniquely determined by
$\sum_{a\in\widetilde{\a}}\mu(a)=1$, whence
 $t=\tfrac1\b$ and $m=\mu$. \hfill\qed\enddemo
\noindent
Note that the proof of proposition 7.1 only used the existence of one cylinder with positive, finite $\mu$-measure.

\bigskip\par
As before, define $F^\natural:X_{\b}\to\Bbb Z_0^{[\b]}=\Z^{[\b]}$ by
$F^\natural(x)_j:=\d_{j,x_1}$.  (Here and throughout $[\cdot]$ is the largest integer {\em lower} bound, i.e. $[\b]<\b$.)
Write $x_1^N\bowtie y_1^N$ if $F^\natural_N(x)=
F^\natural_N(y)$ and $\exists\ k,\ell,\ [x_1^N]\in\widetilde\a_k,\
[y_1^N]\in\widetilde\a_\ell$.

\proclaim{Lemma 7.2}
$\Cal E(X_{\b,0},T,\a)=\{(x,y)\in (X_{\b,0})^2:x_{N+1}^\infty=y_{N+1}^\infty\text{ and }x_1^N\bowtie y_1^N\text{ for}$ $ \text{some }N\}$.
\endproclaim

\heading
\S8 Conformal measures: existence and ergodicity
\endheading
\subheading{Restricted conformal measures}
Fix
$J\subset\{0,1,\dots,[\b]\},\ |J|>1$.  Note that
$X_{\b,0}(J):=X_{\b,0}\cap J^\Bbb N\ne\emptyset$, because
$j_1<j_2\in J\Rightarrow(j_1,j_1,\dots)\in X_{\b,0}(J)$. \par
Setting $\a_n(J):=\{[\u a]\in\a_n:\ \u a\in J^n\}$ and
$\widetilde\a(J):=\bigcup_{n=1}^\infty\widetilde\a\cap\a_n(J)$, we
have that  $(X_{\b,0}(J),S,\widetilde\a(J))$ is a full fibred
system, whence $X_{\b,0}(J)$ is  either a singleton, or
uncountable. In the  latter case $X_{\b}(J):=X_{\b}\cap J^\Bbb
N=\o{X_{\b,0}(J)}$.\footnote{The singleton case is
possible: If $\b=\tfrac{3+\sqrt 5}2$, then $\om=(2,\o{1})$ and
$X_{\b,0}(\{1,2\})=\{(1,1,\dots)\}$. }

As before, we have
$$\Cal E(X_{\b,0}(J),S,\widetilde\a(J)) \subset\goth
T(T|_{X_{\b,0}(J)}) \subset\goth
T(S|_{X_{\b,0}(J)})\subset\goth
G(T|_{X_{\b,0}(J)}).\tag{8.1}$$ Thus, in the
infinite case,  $\goth G(T)\cap [X_{\b,0}(J)\x X_{\b,0}(J)]$ is minimal: Any equivalence class is dense.
\proclaim{Proposition 8.0} Suppose that $X_{\b}(J)$ is
infinite, and let $H:J\to\Bbb R_+$, then
\roster
\item there is a unique $\l>0$ and
$(\tfrac1{\l H},T)$-conformal $\mu_{J,H}\in\Cal P(X_{\b}(J))$;
\item  $\mu_{J,H}$ is a $(S,\widetilde\a)$-product
measure and  is $\Cal E(X_\b)$-invariant and ergodic.
\endroster
\endproclaim\demo{Proof}  In case $\om_\b$ is eventually periodic, $X_{\b}(J)$
is sofic (see \cite{Bl}) and a continuous, equivariant  image of a
TMS. Existence follows from Ruelle's Perron-Frobenius theorem,
which provides a non-atomic, $(\tfrac1{\l H},T)$-conformal
measure.

\par In the case that $\om_\b$ is not eventually periodic, we prove existence (as in \S2 and \S3 of \cite{W}) as follows:
 Endow $X_\b(J)$ with the
lexicographic order topology, disconnecting it at
$\G:=\bigcup_{n\in\Bbb Z}T^n\{(\o 0),\om_\b\}$ to obtain the compact
metric space
$$
Y=Y_\b(J):=\bigl(X_\b(J)\setminus\G\bigr)\cup\bigl(\G\x\{-,+\}\bigr).
$$
There are continuous maps
$T_Y:Y\to Y$ and $\pi:Y\to X_\b(J)$ defined by
$$(T_Y)|_{X_\b(J)\setminus\G}\equiv T,\ \pi|_{X_\b(J)\setminus\G}\equiv \text{\rm Id},\ \
T_Y(\g,\e):=(T\g,\e),\ \pi(\g,\e):=\g,$$  so that $\pi\circ
T_Y=T\circ\pi$. It follows that  $H\circ\pi:Y\to\Bbb R_+$ is
continuous. By Schauder's fixed point theorem, $\exists\ \nu_Y\in\Cal P(Y)$ $(\tfrac1{\l
H\circ\pi}, T_Y)$-conformal for some $\l>0$.

In order to show that $\nu_Y$ descends to a conformal measure on $X_\b(J)$, it is sufficient to show that $\nu_Y(\pi^{-1}\Gamma)=0$, because in this case  $\nu_Y$ is supported
on $X_\b(J)\setminus\Gamma$ so $\pi:(Y,\Cal B(Y),\nu_Y)\to
(X_\b(J),\Cal B,\nu_Y\circ\pi^{-1})$ is a measure theoretic
isomorphism. The measure $\mu_{J,H}:=\nu_Y\circ\pi^{-1}$ is
then the required conformal measure on $X_{\b,0}$.

\medskip
\noindent
{\it Step 1.\/} $\nu_Y\{((\o 0),\pm)\}=0$.
\medskip

\par Otherwise one of
the preimages $z=((j,\o 0),\pm)$ has positive measure.
 The exchangeable orbit of $((j,\o 0),\pm)$ is
$\{((\undersetbrace{n-\text{\rm times}}\to{0,\dots,0},j,\o
0),\pm)\}_{n=1}^\infty$. Since $\nu_Y$ is $(\tfrac1{\l H\circ\pi},
T_Y)$-conformal,  $\nu_Y(\{((\undersetbrace{n-\text{\rm
times}}\to{0,\dots,0},j,\o 0),\pm)\})=\nu_Y(\{z\})\ \forall\ n\ge
0$, a contradiction to the finiteness of $\nu_Y$.

\medskip
\noindent
{\it Step 2.\/} $\nu_Y(\{(\om_\b,\pm)\})=0$.
\medskip

\par Suppose
(to get a contradiction) that $\nu_Y(\{(\om_\b,\pm)\})>0$. In this
case, $\om_\b\in X_\b(J)$, whence $a_n\in J\ \forall\ n\ge 1$. In
particular $a_{min}:=\min\{a_n:n\in\Bbb N\}\in J$. We  construct a
preimage of $(\om_\b,\pm)\})$, which on one hand belongs to
$Y_\b(J)$ (and  by $T_Y$-non-singularity has positive measure),
and on the other hand has an infinite exchangeable orbit with
infinite measure. This is  a contradiction to the finiteness of $\nu_Y$.

Set $ z^{\pm}:=\left((a_{min},\om_\b),\pm\right)$, then $z^\pm\in
Y(J)$  and its exchangeable orbit consists of
$\{(a_{1}^{n-1},a_{min},a_{n+1}^{\infty}):n\in\Bbb N\}$.
This is an infinite set whenever $\om_\b$ is not eventually
periodic.

\medskip
\noindent
{\it Step 3.\/} $\nu_Y(\pi^{-1}\Gamma)=0$.
\medskip
\par For every $\g\in
\Gamma,\ \e=\pm 1$, either $(\g,\e)\not\in Y(J)$ in which case
$\nu_Y\{(\g,\e)\}=0$, or $(\g,\e)\in Y(J)$ and $\exists n,k>0$
such that $T_Y^n((\g,\e))\in T^k\{(\o 0),\om_\b\}$ and
$\nu_Y\{(\g,\e)\}=0$ by $T_Y$--non-singularity of $\nu_Y$.

\medskip
As explained above, having proved that $\nu_Y(\pi^{-1}\Gamma)=0$, we can now obtain a non-atomic, $(\tfrac1{\l H},T)$-conformal
$\mu_{J,H}\in\Cal P(X_{\b,0}(J))$.

\medskip
\par We turn to the  uniqueness part of (1). By $(\tfrac1{\l H},T)$-conformality,
$$\mu_{J,H}(a):=\l^n\widetilde H(a)\ \text{ for all } a\in\widetilde\a_k\cap\a_n,\ \ k,n\ge
1\tag{8.2} $$ where $\widetilde H:\bigcup_{n=1}^\infty\a_n\to\Bbb R_+$ is
defined by
$\widetilde H([a_1^n]):=\prod_{k=1}^nH(a_k)\ \ \ \
([a_1^n]\in\widetilde\a\cap\a_n).$
Noting that if
$B_n:=\sum_{a\in\widetilde\a\cap\a_n}\widetilde H(a)$, then
$\sum_{n=1}^\infty B_n\l^n=1$, we see that the $\l>0$
appearing in (1) is unique.  Equation (8.2) thus determines $\mu_{J,H}$ uniquely. Along the way we have also shown that $\mu_{J,H}$ is
a $(S,\widetilde\a)$-product measure.

Since $\mu_{J,H}$ is a $(S,\widetilde\a)$-product measure, it is
$\Cal E(X_{\b,0},S,\widetilde\a)$-invariant, and (by theorem 2.1,
part 1), $\Cal E(X_{\b,0},S,\widetilde\a)$-ergodic, whence by
(7.3), $\Cal E(X_{\b,0},T,\a)$-ergodic.

To see that it is also $\Cal E(X_{\b,0},T,\a)$-invariance, use
 lemma 7.2 to observe that $\Cal E(X_{\b,0},T,\a)$ is generated by holonomies
of form $(\tilde x,z)\mapsto (\tilde y,z)$ where $[\tilde
x]\bowtie [\tilde y]$. These are all measure preserving.

Finally observe that since $\mu_{J,H}$ is supported on $X_{\b,0}$ and $X_{\b,0}$ is $\Cal E(X_\b,T,a)$ invariant, $\Cal E(X_{\b,0},T,\a)$--ergodicity and invariance is the same as  $\Cal E(X_\b,T,\a)$--ergodicity and invariance.
\hfill\qed\enddemo

\noindent
{\em Remark\/}: The Lebesgue measure $m$
corresponds to $J=\{0,1,\dots,[\b]\}$ , $H=$const.

\proclaim{Corollary 8.1}
 If $X_{\b}(J)$ is infinite, then $\Cal E(X_{\b,0},T,\a)\cap X_{\b}(J)^2$
is topologically transitive.
\endproclaim
\heading{\S 9 From exchangeable measures to conformal measures ($\b$--expansions)}
\endheading
The aim of this section is to prove:
\proclaim{Theorem 9.0}
 If  $\nu \in\goth M(X_{\b,0})$ is locally finite and  $\Cal
E(X_{\b})$-invariant and ergodic, then $\nu=\mu_{J,H}$ for
some $J\subset\{0,1,\dots,[\b]\},\ H:J\to\Bbb R_+$.\endproclaim
\noindent
Theorem 9.0  is false if $\goth M(X_{\b,0})$ is replaced by
$\goth M(X_{\b})$: If
$\b=\frac{1+\sqrt{5}}{2},\ \om_\b=(\o{1,0})$ and $\d_{\om_\b}$ is
$\Cal E(X_{\b})$-invariant and ergodic but not of form
$\mu_{J,H}$ .

\proclaim{Lemma 9.1} Suppose that $(X,T,\a)$ is a fibred
system s.t. $\forall\ x\in X,\ \exists\ N\ge 1$
s.t. $T^N[x_1,\dots,x_N]=X$. If $\psi(x):=\min\,\{N\ge 1:\
T^N[x_1,\dots,x_N]=X\}$ and $Sx:=T^{\psi(x)}x$, then
$\goth T(T)=\goth G(S,\psi)$ (see (1.0)).
\endproclaim\demo{Proof}
\f $\supseteq$:\ \  Suppose that $x\overset{\goth
G(S,\psi)}\to\sim y$, then $\exists\ k,\ \ell\ge 0,\ S^kx=S^\ell
y,\ \psi_k(x)=\psi_\ell(y)=:N$ whence $T^Nx=S^kx=S^\ell y=T^Ny$
and $x\overset{\goth T(T)}\to\sim y$. \f $\subseteq$:\ \ \
 Suppose that $x\overset{\goth T(T)}\to\sim y$, then $\exists\
N\ge 0,\ T^Nx=T^Ny=:z$ and $\exists\ \ell\ge 1$ so that
$$\kappa:=\psi_\ell(z)>N+\max\{\psi(T^jx),\ \psi(T^jy):\ 0\le j\le N\}.$$
It suffices to show that $\kappa=\psi_p(x)=\psi_q(y)$ for some
$p,q\ge 1$ as in this case, $S^px=S^qy$ and $x\overset{\goth
G(S,\psi)}\to\sim y$. \par To see that
$\kappa=\psi_p(x)=\psi_q(y)$ for some $p,q\ge 1$ we prove that
$$T^{N+\kappa}[x_1,\dots,x_{N+\kappa}]=T^{N+\kappa}[y_1,\dots,y_{N+\kappa}]=X.$$
We'll show only that $T^{N+\kappa}[x_1,\dots,x_{N+\kappa}]=X$ (the
other case being analogous). Suppose otherwise, and let
$K:=\max\{k\ge 1:\ T^{N+\kappa-k}[x_{k+1},\dots,x_{N+\kappa}]\ne
X\}$.\par By choice of $\kappa$,
$T^\kappa[x_{N+1},\dots,x_{N+\kappa}]=X$, whence
$T^{\kappa-j}[x_{N+j+1},\dots,x_{N+\kappa}]=X\ \forall\ 0\le j\le
\kappa$. It follows that $K\le N-1$.\par Let $L:=\psi(T^Kx)$, then
$T^L[x_{K+1},\dots,x_{K+L}]=X$. By choice of $\kappa,\
K+L<N+\kappa$ whence
$$\align T^{N+\kappa-K}&[x_{K+1},\dots,x_{N+\kappa}]=\\ &=
T^{N+\kappa-K}([x_{K+1},\dots,x_{K+L}]\cap
T^{-L}[x_{K+L+1},\dots,x_{N+\kappa}])\\ &=
 T^L[x_{K+1},\dots,x_{K+L}]\cap T^{N+\kappa-(K+L)}[x_{K+L+1},\dots,x_{N+\kappa}]\\
& =T^{N+\kappa-(K+L)}[x_{K+L+1},\dots,x_{N+\kappa}]\ne
X\endalign$$ contradicting maximality of $K$.
 \hfill\qed\enddemo
\proclaim{Lemma 9.2}
$\Cal E(X_{\b,0},T,\a)=\goth
G(S|_{X_{\b,0}},\Phi^{(\b)})$ where
$\Phi^{(\b)}:=(F^\flat,\psi):X_{\b,0}\to\Bbb G_\b:=\Bbb
Z^{[\b]}\x\Bbb Z$ and
$F^\flat:=\sum_{k=0}^{\psi-1}F^\natural\circ T^k$.
\endproclaim\demo{Proof} By proposition 2.0,
$(x,y)\in\Cal E(X_{\b,0},T,\a)$ iff $(x,y)\in\goth T(T)\ \&\
F^\natural_n(x)=F^\natural_n(y)$ whenever $T^nx=T^ny.$
By lemma 9.1, $(x,y)\in\goth T(T)$ iff $\exists\ k,\ell\ge 1$ with
 $\psi_k(x)=\psi_\ell(y)=:N$ and $S^kx=T^Nx=T^Ny=S^ky$. Thus
$(x,y)\in\Cal E(X_{\b,0},T,\a)$ iff $\exists\ k,\ell\ge
1$ such that
\roster
\item $\psi_k(x)=\psi_\ell(y)=:N$,
\item $S^kx=T^Nx=T^Ny=S^ky$, and
\item $F^\flat_k(x)=
F^\natural_N(x)=F^\natural_N(y)=F^\flat_\ell(y)$.
\endroster
Equivalently
$(x,y)\in\goth G(S,\Phi^{(\b)})$. \hfill\qed\enddemo
\subheading{Proof of theorem 9.0}
By (7.3), $\nu$ is $\Cal
E(X_{\b,0},S,\widetilde\a)$-invariant and $\goth
T(S|_{X_{\b,0}})$-ergodic. By proposition 2.2, $\nu$ is recurrent, and proportional to a $(S,\widetilde\a)$-product measure. In particular:
\roster
\item for
each $0\le b\le [\b]$, $N_b:=\sum_{k=1}^\infty\d_{x_k,b}=0,\infty$
$\nu$-a.e.;
\item $\nu$ is either a point mass or non-atomic;
\item $\nu$ is $\goth G(S|_{X_{\b,\nu}})$-non-singular where $X_{\b,\nu}:=$ the
closed support of $\nu$ in $X_{\b}$.
\endroster

Let $J:=\{b\in\{0,1,\dots,[\b]\}:\ N_b=\infty\}$.  Either $X_\b(J)$ is a singleton and $\nu$ is a point mass, or $X_\b(J)$ is uncountable and  $\Cal E(X_{\b,0},T,\a)\cap X_{\b}(J)^2$ is topologically transitive (corollary 8.1). The first case is covered by the theorem. In the second case  $X_{\b,\nu}=X_\b(J)$, so $\nu$ is not a point mass, whence by (2) non-atomic. Henceforth we restrict ourselves to this case.

\par Fix $j_0:=\min J$, and  define
$$
\align
F^{J,\natural}:X_\b(J)\to\Bbb Z^{J\setminus\{j_0\}} &,\text{ where }  F^{J,\natural}(x)_j:=\d_{x_1,j}\\
F^{J,\flat}:X_\b(J)\to\Bbb
Z^{J\setminus\{j_0\}} &, \text{ where } F^{J,\flat}:=\sum_{k=0}^{\psi-1}F^{J,\natural}\circ T^k
\endalign
$$
$$
\Phi=\Phi^{(J)}:X_\b(J)\to\Bbb G\textrm{ where }\Bbb G:= \Z^{J\setminus\{j_0\}}\x\Bbb Z \textrm{ and }
\Phi:=(F^{J,\flat},\psi).
$$
By lemma 9.2
$ \Cal E(X_\b(J),T,\a) =\goth G(S|_{X_\b(J)},\Phi)\cong
\goth G(S_{\Phi}|_{X_\b(J)\x\Bbb
G})\cap(X_\b(J)\x\{0\})^2.$ By proposition 1.0,
there is a unique $\s$-finite, $\goth G(S_{\Phi}|_{X_\b(J)\x\Bbb
G})$-invariant, ergodic measure $m$ on $X_\b(J)\x\Bbb G$
so that
$$m(A\x\{0\})=\nu(A)\ \text{ for all } A\in\B(X_\b(J)).$$

\medskip
\noindent
{\it Step 1.\/} $m(X_\b(J)\x\{g\})<\infty\ \forall\ \ g\in\Bbb G$.
\medskip

We claim first that $\nu_g\ll\nu\ \forall\ g\in\Bbb G$ where
$\nu_g(A):=m(A\x\{g\})$. To see this, suppose that
$B\in\B(X_\b(J)),\ m(B\x\{g\})>0$, then by $\goth G(S_\Phi|_{X_\b(J)\x\Bbb
G})$-ergodicity of $m,\ \exists\ k,\ell\ge 1,\
a\in\widetilde\a_k,\ b\in\widetilde\a_\ell$ and $A\in\B(X_\b),\
\nu(A)>0,\ A\subset a$ so that $\Pi(A\x\{0\})\subset B\x\{g\}$
where $\Pi:a\x\Bbb G\to b\x\Bbb G$ is defined by
$$\Pi(y,z):=(\pi(y),z+\widehat{\Phi}(y,\pi(y))),\ \ \pi(a,x):=(b,x).$$
In particular $\Pi(A)\subset B$ whence, by $\goth
G(S|_{X_\b(J)})$-non-singularity of $\nu,\ \nu(B)>0$.
\par Next, for $g\in\Bbb G$, define $Q_g:X_\b(J)\x\Bbb G\to X_\b(J)\x\Bbb G$ by $Q_g(x,z):=(x,z+g)$.
Evidently $(Q_g\x Q_g)[\goth G(S_\Phi|_{X_\b(J)\x\Bbb G})]=\goth
G(S_\Phi|_{X_\b(J)\x\Bbb G})$ whence $m\circ Q_g$ is also a
$\s$-finite, $\goth G(S_\Phi|_{X_\b(J)\x\Bbb G})$-invariant,
ergodic measure $m$ on $X_\b(J)\x\Bbb G$. It follows that either
$m\circ Q_g\perp m$ or $m\circ Q_g=c_gm$ for some $c_g\in\Bbb
R_+$. In particular,
$$m(X_\b\x\{g\})=\cases & 0\ \ \ \ \ \ \ m\circ
Q_g\perp m,\\ & c_g<\infty\ \ \ \ \ \ \text{\rm else.}\endcases$$

\medskip
\noindent
{\it Step 2.\/} $\Bbb H:=\{h\in\Bbb
G:\ m\circ Q_h\sim m\}$ is a subgroup, and $m\circ
Q_h\perp m\ \forall\ h\notin \Bbb H$. We claim that it is enough to show that $\Bbb H=\Bbb G$.
\medskip
Indeed, if $\Bbb H=\Bbb G$, then
$m(A\x\{g\})=e^{H(g)}\mu(A)$ for some homomorphism $H:\Bbb
G\to\Bbb R$, whence $\tfrac{d\mu\circ S}{d\mu}=e^{-H\circ\Phi }$
and $\mu=\mu_{J,h}$ where $h:=H|_{\Bbb
Z^{J\setminus\{j_0\}}\x\{0\}}$.

\medskip
\noindent
{\it Step 3.\/} $\<\Phi(X_\b(J))\>=\Bbb G$ and   $\Phi(X_{\b}(J))\subset \Bbb H$. Consequently, $\Bbb H=\Bbb G$
\medskip
Let $\{e_i\}$ be the standard basis for $\Z^{J\setminus\{j_0\}}$.
To see the first identity, note that for every
$j\in J,\ j\ne [\b]$, $[j]\in\widetilde\a$ and so
$(e_j,1)=\Phi([j])\in\Phi(X_\b(J))\ \ (j\in
J\setminus\{j_0,[\b]\})$ and $(0,1)=\Phi([j_0])\in\Phi(X_\b(J))$.
\par In case $[\b]\in J$,
$\#\{n\ge 1:\ \om_n=[\b]\}=\infty$ and $\exists\ N\ge 1$ with $
\om_1^N=([\b],\om_2,\dots, \om_{N-1},[\b])$ where $\om_2,\dots,
\om_{N-1}\ne [\b]$, whence $[[\b],j_0^k]\in\widetilde\a$ for some
$1\le k\le N-1$ and
$(e_{[\b]},k+1)=\Phi([[\b],j_0^k])\in\Phi(X_\b(J))$.

Next, we show  that $\Phi(X_{\b}(J))\subset \Bbb H$. Accordingly,
 let $\Phi\equiv h$ on $a\in\widetilde\a$, and assume by way of contradiction
that  $h\notin \Bbb H$. In this case,
 $m\circ Q_h\perp m$ and $\exists\
Z\in \B(X),\ m(Z\x\{0\})=1,\ m(Z\x\{-h\})=0$. It follows that
 $\exists\ \xi\subset\bigcup_{n=1}^\infty\widetilde\a_n$
countable, such that $U:=\bigcup_{A\in\xi}A\supset Z$ and
$m(U\x\{-h\})<\tfrac13$.

Since
$m(U\x\{0\})=1,\ \exists\ A\in\xi$ such that
$m(A\x\{-h\})<\tfrac13m(A\x\{0\})$. \bigskip\par Define
$\kappa:A\to A$ by $\kappa(A,x)\ :=\ (A,a,x)$. Evidently
$(x,\kappa(x))\in\goth G(S)\ \forall\ x\in A$, and
 $\widetilde\Phi(x,\kappa(x))=-h\ \forall\ x\in A$ . Thus, if
$\widetilde\kappa:A\x\Bbb G\to A\x\Bbb G$ is defined by
$\widetilde\kappa(x,g):=(\kappa(x),g-h)$, then
$((x,g),\widetilde\kappa(x,g))\in\goth G(S_{\Phi})\ \forall
(x,g)\in A\x\Bbb G$. It follows that
$$m(A\x\{0\})=m(\widetilde\kappa(A\x\{0\}))\le m(A\x\{-h\})<\tfrac13m(A\x\{0\}).$$
By the contradiction, $h\in \Bbb H$, proving the step and with it the theorem.
\qed

\bigskip

\centerline{\bf References}

\medskip
\noindent [ADU]  J. Aaronson, M. Denker, M. Urba\'nski,  Ergodic theory for Markov fibred systems and parabolic rational maps.  {\it Trans. Amer. Math. Soc.}  {\bf 337}  (1993),  no. 2, 495--548.

\medskip
\noindent [AD] J. Aaronson, M.  Denker, Group extensions of
Gibbs-Markov maps,   {\it Probab. Theory Related Fields}  {\bf
123} (2002), no. 1, 38--40.

\medskip
\noindent [ADSZ] J. Aaronson, M. Denker, O.  Sarig, R.
Zweim\"uller, Aperiodicity of cocycles and conditional local limit
theorems,  {\it Stoch. Dyn.}  {\bf 4}  (2004),  no. 1, 31--62.

\medskip
\noindent [ANSS] J. Aaronson, H.  Nakada, O.  Sarig, R.  Solomyak,
 Invariant measures and asymptotics for some skew products.
{\it Israel J. Math.} {\bf 128} (2002), 93--134.

\medskip
\noindent [ANSS1] ------------- , Corrections to :Invariant
measures and asymptotics for some skew products,  {\it Israel J.
Math.} {\bf 138} (2003), 377--379.

\medskip
\noindent [Bl]\ \  F. Blanchard, $\beta$-expansions and symbolic
dynamics, {\it Theoret. Comput. Sci.} {\bf 65},  (1989), no. 2,
131--141.

\medskip
\noindent [Ch]\ \  K.L Chung, {\it Markov chains with stationary
transition probabilities}, Springer, Heidelberg 1960.

\medskip
\noindent [D-F]\ \ P. Diaconis, D. Freedman, De Finetti's theorem
for Markov chains,  {\it Annals Probab.} {\bf 8}, no. 8, (1980)
115-130.

\medskip
\noindent [E]\ \  E. Effros,  Transformation groups and $C\sp{*}
$-algebras, {\it Ann. of Math.} (2) {\bf 81}, (1965) 38--55.

\medskip
\noindent [F-M]\ \ J. Feldman and C. C. Moore, Ergodic equivalence
relations, cohomology, and von Neumann algebras.I {\it Trans. Am.
Math. Soc. }, Volume 234, {\bf 2}, (1977), 289--324.

\medskip
\noindent [Fu]\ \  L.  Fuchs,  {\it Abelian groups.} International
Series of Monographs on Pure and Applied Mathematics Pergamon
Press, New York-Oxford-London-Paris 1960.

\medskip
\noindent [Ge]\ \ A. O.   Gel'fond, A common property of number
systems. (Russian){ Izv. Akad. Nauk SSSR. Ser. Mat. }{\bf 23},
(1959) 809--814.

\medskip
\noindent [Gl]\ \  J. Glimm,  Locally compact transformation
groups, {\it Trans. Amer. Math. Soc.} {\bf 101}, (1961) 124--138.

\medskip
\noindent [G-S]\ \ G. Greschonig, K. Schmidt, Ergodic decomposition of quasi-invariant probability measures,
{\it Colloq. Math.} {\bf  84/85}  (2000), part 2, 495--514.

\medskip
\noindent [Gr]\ \   L.A. Grigorenko, On the $\sigma$-algebra of
symmetric events for a countable Markov chain.  {\it
Theory Probab. Appl.}{\bf 24} (1979) 199-204.

\medskip
\noindent [H-S]\ \  E. Hewitt, L.J.  Savage, Symmetric measures on
Cartesian products, {\it Trans. Amer. Math. Soc.} {\bf 80} (1955),
470--501.

\medskip
\noindent [I-Ta]\ \ S. Ito, Y.  Takahashi, Markov subshifts and
realization of $\beta $-expansions, {\it J. Math. Soc. Japan},
{\bf 26}, (1974), 33--55.

\medskip
\noindent [L]\ \  A. N. Liv\v{s}ic, Certain properties of the
homology of $Y$-systems. {\it Mat. Zametki} {\bf  10}
(1971), 555--564. Engl. Transl. in {\it Math. Notes} {\bf 10}
(1971), 758--763.

\medskip
\noindent [Ma]\ \ D. Maharam, Incompressible transformations, {\it
Fund. Math.}  {\bf 56}, (1964), 35-50.

\medskip
\noindent [Me]\ \  P-A. Meyer, {\it Probability and potentials,}
Blaisdell Publishing Co. Ginn and Co., Waltham, Mass.-Toronto,
Ont.-London 1966.

\medskip
\noindent [Pa]\ \  W. Parry, On the $\beta $-expansions of real
numbers,{\it Acta Math. Acad. Sci. Hungar.}{\bf 11} (1960)
401--416.

\medskip
\noindent [Pa-S]\ \  W. Parry, K. Schmidt,  Natural coefficients
and invariants for Markov-shifts, {\it Invent. Math.} {\bf 76}
(1984), no. 1, 15--32.

\medskip
\noindent [Pe-S]\ \  K. Petersen, K. Schmidt, Symmetric Gibbs
measures. {\it Trans. Amer. Math. Soc.} {\bf 349}, (1997),
2775--2811.

\medskip
\noindent [Re]\ \ A. R\'{e}nyi, Representations for real numbers
and their ergodic properties,{\it Acta Math. Acad. Sci.
Hungar.}{\bf 8} (1957) 477--493.

\medskip
\noindent [Sa]\ \  O. Sarig, Existence of Gibbs measures for
countable Markov shifts, {\it Proc. Amer. Math. Soc.}{\bf 131}
(2003), no. 6, 1751--1758

\medskip
\noindent [Schm]\ \ K. Schmidt,  Infinite invariant measures on
the circle. {\it Symposia Mathematica}, Vol. {\bf XXI} (Convegno
sulle Misure su Gruppi e su Spazi Vettoriali, Convegno sui Gruppi
e Anelli Ordinati, INDAM, Rome, 1975), pp. 37--43. Academic Press,
London, 1977.

\medskip
\noindent [Schw]\ \  F. Schweiger, {\it Ergodic theory of fibred
systems and metric number theory}, Oxford Science Publications.
The Clarendon Press, Oxford University Press, New York, 1995.

\medskip
\noindent [Sm]\ \ M. Smorodinsky, $\b$--automorphisms are
Bernoulli shifts, {\it Acta. Math. Acad. Sci. Hungar.} {\bf 24}
(1973), 273--278.

\medskip
\noindent [VJ]\ \ D. Vere-Jones, Ergodic properties of nonnegative
matrices. I. {\it Pacific J. Math.} {\bf 22} 1967 361--386.

\medskip
\noindent [W]\ \ P. Walters,  Equilibrium states for $\beta
$-transformations and related transformations. {\it Math. Z.} {\bf
159} (1978), no. 1, 65--88.

\enddocument

{\  \medskip
\noindent [L-M]\ \ A. Louveau, G. Mokobodzki,  On measures ergodic
with respect to an analytic equivalence relation, {\it Trans.
Amer. Math. Soc.}{\bf 349} (1997), no. 12, 4815--4823.
\medskip
\noindent [C] I.V.   Culanovskij, On cycles in Markov chains
(English) {\it Select. Transl. Math. Stat. Probab.}{\bf 1} (1961)
1-5.

\medskip
\noindent [Gu]\ \  B.M. Gurevich, The local limit theorem for
Markov chains and regularity conditions (Russian, English) {\it
Theor. Probab. Appl.}{\bf 13}(1968) 182-188 ; translation from
{\it Teor. Veroyatn. Primen.}{\bf 13} (1968) 183-190.
\medskip
\noindent [H-S]\ \ A.N.  Kolmogorov, A local limit theorem for
Markov chains (English) {\it Select. Transl. Math. Stat. Probab.
}{\bf 2} (1962) 109-129.
\medskip
\noindent [Ro]\ \ I.Z. Rozenknop, On some properties of the
totality of closed paths in a system of n states and given
transitions among them (English) {\it Select. Transl. Math. Stat.
Probab. }{\bf 1} (1961) 7-12.
\par We
show first that $\exists$ a $\mu$- sweep-out set of finite
$\mu$-measure. To this end, define $N_a:=\sum_{k=0}^\infty
1_a\circ T^k\ \ (a\in\a)$. We claim that $\mu([1\le
N_a<\infty])=0\ \forall\ a\in\a$ since otherwise, for some
$s\in\a,\ b\in\b,\ k\in\Bbb N,\ \exists\ K\subset\Bbb N,\ |K|=k$
such that $\mu(A_{K,s})=:c>0$ where
$$A_{K,s}:=\{x=(a_1,a_2,\dots)\in b:\ a_n=s\ \forall\ n\in K,\ a_n\ne
s\ \forall\ n\notin K\}.$$ Evidently for each $L\subset\Bbb N,\
|L|=k,\ A_{L,s}\overset{\Cal E(X,T,\a)}\to\rightarrow A_{K,s}$ and
$\mu(A_{L,s})=c$. This contradicts so By ,  Fix $b\in\b$ and let
$b\in\a_n$, then $a\overset{\Cal E(X,T,\a)}\to\rightarrow b\
\forall\ a\in\a_n$ and $\mu$ satisfies the precondition of
proposition 1.0. Accordingly we assume (without loss of
generality) that $\mu$ is (in addition) $\Cal E(X,T,\a)$-ergodic.

\par In
many interesting cases, we have $\o{\pi(X)}\setminus\pi(X)$
countable. For example, if $X=[0,1),\ Tx=2x\ \mod 1$ and
$\a=\{a_0,a_1\}=\{[0,\tfrac12),[\tfrac12,1)\}$ then $(X,T,\a)$ is
a   fibred system with $\Sigma(X,T,\a)=\{0,1\}^\Bbb N$ (the full
shift) and $\o{\pi(X)}\setminus\pi(X)=\{x\in\Sigma(X,T,\a):\
x\rightarrow 1\}.$ The condition is also necessary, see
\cite{??}.\demo{Proof}\par Define $\Cal A:\bigcup_{n=2}^\infty
S^n\to \{0,1\}$ by $$\Cal
A(s_1,\dots,s_n):=\prod_{k=1}^{n-1}A_{s_k,s_{k+1}},$$ then
$[s_1,\dots,s_n]\ne\emptyset$ in $X$ iff $\Cal
A(s_1,\dots,s_n)=1$.\par Fix $s\in S$ and let
$$\widetilde\a_{s,n}:=\{a=[s,a_1,\dots,a_{n-1},s]:\
a_1,\dots,a_{n-1}\in S\setminus\{s\},\ \Cal A(a)=1\},\
\widetilde\a_{s}:=\bigcup_{n=1}^\infty\widetilde\a_{s,n}.$$ Define
$H:\widetilde\a_{s}\to\Bbb R_+$ by
$$H([s,a_1,\dots,a_{n-1},s]):=h(s)\prod_{k=1}^{n-1}h(a_k).$$
Evidently, for $n\ge 1$,
$$B_n:=\sum_{a\in\widetilde\a_{s,n}}H(a)\le M^n$$
and it follows that the power series $\sum_{n=1}^\infty B_nz^n$
has a positive radius of convergence. Thus $\exists\ !\ c>0$ so
that $\sum_{n=1}^\infty B_nc^n=1$. Set $\o h:=ch$ and define
$p:\widetilde\a_{s}\to\Bbb R_+$ by
$$p([s,a_1,\dots,a_{n-1},s]):=\o h(s)\prod_{k=1}^{n-1}\o h(a_k)
=c^nH([s,a_1,\dots,a_{n-1},s]).$$ Evidently, $p\in\Cal
P(\widetilde\a_{s})$ as
$\sum_{a\in\widetilde\a_{s}}p(a)=\sum_{n=1}^\infty B_nc^n=1$.
\par Define $\nu\in\Cal P([s])$ as the
$([s],T_{[s]},\widetilde\a_{s})$-product measure with marginal
$p$. We claim that $\nu$ is $\Cal E(X,T,\a)\cap [s]\x
[s]$-invariant, ergodic, and also has the property that
$$\sum_{n=0}^\infty 1_{[a]}\circ T^n=\infty\ \nu-\text{\rm a.e. on
}[s]\ \forall\ a\in \bigcup_{n=1}^\infty S^n,\ \Cal A(a)=1.\tag1$$
The $\Cal E(X,T,\a)\cap [s]\x [s]$-invariance of $\nu$ follows
from its product form.
\par By theorem 2.1 (part 1), $\nu$ is $\Cal
  E([s],T_{[s]},\widetilde\a_s)$-ergodic, whence $\Cal
E(X,T,\a)\cap [s]\x [s]$-ergodic since $\Cal
E([s],T_{[s]},\widetilde\a_s)
  \subset\Cal E(X,T,\a)\cap [s]\x [s]$. \par To see (1), let $\kappa\ge
1,\ a\in  S^\kappa,\ \Cal A(a)=1$. Suppose that $\#\{1\le k\le
\kappa:\ a_k=s\}=N$, then $\exists\
b_1,b_2,\dots,b_{N+1}\in\widetilde\a_s^N$ so that if
$$[b_1,b_2,\dots,b_{N+1}]=[c_1,c_2,\dots,c_M]\in\a_M$$
then  for some $k,\ a_i=c_{k+i}\ \ (1\le i\le \kappa)$. \par By
the ergodic theorem for $([s],\B([s]),\nu,T_{[s]}),\ \nu$ a,e, on
$[s]$:

$$\sum_{j=0}^\infty 1_{[b_1,b_2,\dots,b_{N+1}]}\circ
T_{[s]}^j=\infty,\ \ \text{\rm  whence }\sum_{n=0}^\infty
1_{[a]}\circ T^n=\infty.$$
\par Now let $\mu\in\goth M(X)$ be the unique $\Cal E(X,T,\a)$-invariant, ergodic
 measure satisfying $\mu|_{[s]}\equiv\nu$. The required properties (Markovianity and conformality)
 will be established by showing $\mu\in\goth M_\a(X)$.
 \par To see this, note that by topological transitivity of $\Cal E(X,T,\a)$,
  $\mu$ is globally supported and locally finite. Also,
$$\sum_{n=0}^\infty 1_{[a]}\circ T^n=\infty\
\mu-\text{\rm a.e.}\ \forall\ a\in \bigcup_{n=1}^\infty S^n,\ \Cal
A(a)=1\tag2$$ since this property is $\Cal E(X,T,\a)$-invariant
(even $\goth G(T)$-invariant).
  \par If $t\in S$, then
  $\mu|_{[t]}$ is

  \smallskip\f$\bullet\ \ \ \ $ locally finite;
  \smallskip\f$\bullet\ \ \ \ $
  $\Cal E(X,T,\a)\cap [t]\x [t]$-invariant, ergodic;
  \smallskip\f $\bullet\ \ \ \ $
  recurrent with respect to $([t],T_{[t]},\widetilde\a_t)$ (by (2)). \par Now
  $$\Cal E([t],T_{[t]},\widetilde\a_t)
  \subset\Cal E(X,T,\a)\cap [t]\x [t]\subset\goth T(T_{[t]}),$$ and so $\mu|_{[t]}$ is
 also $\Cal E([t],T_{[t]},\widetilde\a_t)$-invariant and $\goth T(T_{[t]})$-ergodic.
 \par By $\Cal E([t],T_{[t]},\widetilde\a_t)$-invariance,
 $\mu_{[t]}$ is a mixture of $\Cal E([t],T_{[t]},\widetilde\a_t)$-product measures
 (proposition 2.2).
 By $\goth T(T_{[t]})$-ergodicity, there can be only one such component
 $\Cal E([t],T_{[t]},\widetilde\a_t)$-product measure.
 Thus $\mu([t])<\infty$.\hfill\qed\enddemo

Greschonig, Gernot; Schmidt, Klaus. Ergodic decomposition of quasi-invariant probability measures. Dedicated to the memory of Anzelm Iwanik. Colloq. Math.  84/85  (2000), , part 2, 495--514.}\enddocument